\documentclass[aap,MSNbibl,citesort,dvips]{arximspdf}
\usepackage{subenv}

\doi{10.1214/11-AAP764}
\volume{22}
\issue{1}
\pubyear{2012}
\firstpage{172}
\lastpage{212}

\makeatletter

\newproclaim{assumption}{Assumption}[section]
\newproclaim{example}{Example}[section]
\newproclaim{definition}{Definition}[section]
\newproclaim{remark}{Remark}[section]

\newtheorem{lemma}{Lemma}[section]
\newtheorem{theorem}{Theorem}[section]
\newtheorem{proposition}{Proposition}[section]

\newcommand{\dd}{\,d}

\makeatother

\begin{document}
\begin{frontmatter}

\title{Asymptotics of robust utility maximization}
\runtitle{Asymptotics of robust utility maximization}

\begin{aug}
\author[A]{\fnms{Thomas} \snm{Knispel}\corref{}\ead[label=e1]{knispel@stochastik.uni-hannover.de}\thanksref{t1}}
\runauthor{T. Knispel}
\affiliation{Leibniz Universit\"{a}t Hannover}
\address[A]{Institut f\"{u}r Mathematische Stochastik\\
Leibniz Universit\"{a}t Hannover\\
Welfengarten 1\\
30167 Hannover\\
Germany\\
\printead{e1}} 
\end{aug}

\thankstext{t1}{Supported by Deutsche Forschungsgemeinschaft through
International Research Training Group 1339
``Stochastic Models of Complex Processes'' (Berlin).}

\received{\smonth{3} \syear{2010}}
\revised{\smonth{11} \syear{2010}}

%
\begin{abstract}
For a stochastic factor model we maximize the long-term growth rate of
robust expected power utility with parameter $\lambda\in(0,1)$. Using
duality methods the problem is reformulated as an infinite time
horizon, risk-sensitive control problem. Our results characterize the
optimal growth rate, an optimal long-term trading strategy and an
asymptotic worst-case model in terms of an ergodic Bellman equation.
With these results we propose a duality approach to a ``robust
large deviations'' criterion for optimal long-term
investment.
\end{abstract}

%

\end{frontmatter}

\section{Introduction}
One of the basic tasks in mathematical finance is to choose an
``optimal'' payoff among all available financial
positions which are affordable given an initial capital endowment.
In mathematical terms, a~payoff at a terminal time corresponds to a
real-valued random variable on some measurable space $(\Omega
,\mathcal{F})$ and an investor faces a set $\mathfrak
{X}$ of such financial positions.
Any
formulation of optimality will involve the investor's individual
\textit{preferences} $\succ$ on $\mathfrak
{X}$. The relation $X\succ Y$ means that the investor prefers the
payoff $X$ over $Y$. Under mild conditions such preferences admit a
\textit{numerical representation} $\mathcal{U}\dvtx\mathfrak{X}\rightarrow
\mathbb{R}$ (see, e.g., \cite{foellmerschied}); that is, for $X,Y\in
\mathfrak{X}$ it holds that
\[
X\succ Y \quad\Longleftrightarrow\quad\mathcal{U}(X)>\mathcal{U} (Y).
\]
In this context, Savage \cite{savage} clarified the conditions which
guarantee that a~preference order admits the specific numerical representation
%
\begin{equation}
\label{eq:gl149}
\mathcal{U}(X)=E_Q[u(X)]=\int u(X(\omega))Q(d\omega),\qquad X\in
\mathfrak{X},
\end{equation}
in terms of an increasing continuous function $u\dvtx \mathbb{R}\rightarrow
\mathbb{R}\cup\{-\infty\}$ and a probability measure $Q$ on $(\Omega
,\mathcal{F})$. Here $Q$ appears as a ``subjective''
probability measure which is implicit in the investor's
preferences, and which may differ from a given ``objective''
probability measure. The function $u$ in (\ref{eq:gl149})
will be concave if the investor is assumed to be risk averse. In that
case, $u$ is called a \textit{utility function}.

The literature on optimal investment decisions in a financial market
usually involves the maximization of a utility functional (\ref
{eq:gl149}) with respect to a~given measure~$Q$. Typically, $Q$ is
assumed to model the evolution of future stock prices and is thus
viewed as the ``objective'' measure. But the
price dynamics are not really known accurately, and so the choice of
the evaluation measure $Q$ is itself subject to \textit{model
uncertainty} or \textit{model ambiguity}, also called \textit{Knightian
uncertainty} in the economic literature. There is another reason to
depart from the standard setting of expected utility as formulated in
(\ref{eq:gl149}): some very plausible preferences such as the famous
\textit{Ellsberg paradox} are not consistent with (\ref{eq:gl149}) (see,
e.g., \cite{foellmerschied}, Example~2.75). In order to overcome this
limitation, Gilboa and Schmeidler \cite{gilboaschmeidler} proposed a~more flexible set of axioms for preference orders which leads to a
``robust'' extension of~(\ref{eq:gl149}): instead of
a single measure $Q$ the numerical representation of the preference
order involves a whole class $\mathcal{Q}$ of probability measures and
takes the form of a ``coherent'' robust utility functional
%
\begin{equation}
\label{eq:gl150}
\mathcal{U}(X)=\inf_{Q\in\mathcal{Q}}E_Q[u(X)].
\end{equation}
This representation suggests the following interpretation: the investor
has in mind a collection of possible probability distributions of
market events and takes a worst-case approach in evaluating the
expected utility of a given payoff. In recent years, there is an
increasing interest in the maximization of the robust expected utility
(\ref{eq:gl150}) of wealth $X_T^{x_0,\xi}$ attainable at time $T>0$ by
investing in a financial market using some self-financing trading
strategies $\xi$ and the initial capital $x_0$
%
\begin{equation}
\label{eq:gl534}
\mbox{maximize } \inf_{Q\in\mathcal{Q}}E_Q[u(X_T^{x_0,\xi})]
\mbox{ among all self-financing strategies } \xi.
\end{equation}
For general semimartingale models, this optimization problem can be
solved by a \textit{duality approach} (sometimes also called martingale
approach) (see, e.g., Quenez~\cite{quenez}, Schied and Wu \cite
{schiedwu} or F\"{o}llmer and Gundel \cite{foellmergundel}). Their
results provide a robust extension of the seminal paper by Kramkov and
Schachermayer \cite{kramkovschachermayer} for the classical utility
maximization problem in incomplete markets. The main advantage of the
duality approach lies in the fact that the primal saddle-point problem
is reduced to a minimization problem on the dual side. In many cases,
the dual problem is much simpler and can be tackled with another
optimization technique (dynamic programming, BSDE).

For a finite maturity, however, the optimal investment strategies for
(\ref{eq:gl534}) will typically be time dependent, and they are often
difficult to compute. Instead we propose an asymptotic approach: we
consider a long-term investment model with one riskless and one risky
asset whose drift coefficients are affected by an external \textit
{stochastic factor process} of diffusion type. Our model takes into
account ambiguity about the ``true'' drift terms
of both the factor process and the risky asset. The class $\mathcal{Q}$
of possible prior models corresponds to affine perturbations of the
drift terms in a given reference model and is parameterized by
stochastic controls. In this paper we focus on \textit{power utility}
$u(x)=\frac{1}{\lambda}x^\lambda$ with parameter $\lambda\in(0,1)$,
but other utility functions are also feasible (cf. Remark \ref
{th:negpow}). In our model the robust expected power utility will grow
exponentially as time $T\uparrow\infty$, and this suggests to
%
\begin{equation}
\label{eq:gl532}
\mbox{maximize } \mathop{\overline{\lim}}_{T\uparrow\infty
} \frac{1}{T}\ln\inf_{Q\in\mathcal{Q}}E_{Q}[(X_T^{x_0,\xi
})^\lambda] \mbox{ among all strategies } \xi.
\end{equation}
This asymptotic formulation has the advantage of allowing for
stationary optimal policies and may thus be more tractable. On the
other hand, the asymptotic ansatz provides useful insight for portfolio
management with long but finite time horizon.

For the nonrobust case $\mathcal{Q}=\{Q\}$, problem (\ref{eq:gl532})
is closely related to the maximization of the portfolio's \textit
{risk-sensitized expected growth rate},
%
\begin{equation}
\label{eq:gl354}
\Lambda_Q(\theta,\xi):=\mathop{\underline{\lim}}_{T\uparrow\infty
} -\frac{2}{\theta T}\ln E_Q\biggl[\exp\biggl(-\frac{\theta}{2} \ln
X_T^{x_0,\xi}\biggr)\biggr],\qquad \theta\not=0.
\end{equation}
In order to explain the nature of this criterion, let us consider the
\textit{entropic monetary utility functional} $
\mathcal{U}_\theta(X):=-\frac{2}{\theta}\ln E_{Q}[\exp(-\frac{\theta
}{2}X)]$, where $\theta$ is a positive constant. The functional
$\mathcal{U}_\theta$ is also well defined for $\theta<0$, and it can be
extended to $\theta=0$ via $\mathcal{U}_0(X):=\lim_{\theta\rightarrow
0}\mathcal{U}_\theta(X)=E_{Q}[X]$. A Taylor expansion around $\theta=0$
(cf., e.g., \cite{whittle2}, page 5) yields
%
\begin{equation}
\label{eq:gl555}
\mathcal{U}_\theta(X)=E_{Q}[X]+\frac{\theta}{4}\operatorname
{Var}_{Q}[X]+O(\theta^2).
\end{equation}
Thus $\theta$ can be interpreted as a ``risk sensitivity''
parameter that weights the impact of variance. In
particular, the Taylor expansion (\ref{eq:gl555}) suggests that
\[
\Lambda_Q(\theta,\xi)=\mathop{\underline{\lim}}_{T\uparrow\infty
} \frac{1}{T}E_{Q}[\ln X_T^{x_0,\xi}]+\frac{\theta
}{4}\mathop{\underline{\lim}}_{T\uparrow\infty} \frac
{1}{T}\operatorname{Var}_{Q}[\ln X_T^{x_0,\xi}].
\]
The first term at the right-hand side is the portfolio's \textit
{risk-neutral expected growth rate}. The second term provides a risk
adjustment specified by the portfolio's asymptotic variance and the
risk sensitivity parameter $\theta$, and so~$\Lambda_Q(\theta,\xi)$ can
indeed be seen as the risk-sensitized expected growth rate of wealth.
On the other hand, the long-run growth rates of expected power utility
$u(x)=(\theta/2)x^{\theta/2}$ are, up to constants, of the form $\Lambda
_Q(\theta,\xi)$, and the limit $\theta\rightarrow0$ corresponds to the
growth rate of expected logarithmic utility. Such risk-sensitized
portfolio optimization problems on an infinite time horizon have
received much attention (see, e.g.,
\cite{flemingsheu,flemingsheu1,bieleckipliska1,bieleckipliska2,pham,nagai2,kurodanagai}). In those papers, the
maximization of (\ref{eq:gl354}) among a class of trading strategies,
viewed as dynamic controls, is reformulated as an infinite time
horizon, risk-sensitive control problem of the kind studied in Fleming
and McEneaney \cite{flemingmceneaney}. The rewritten problem leads to
an auxiliary finite horizon ``exponential of integral criterion.''
This is a standard problem in stochastic control theory, and its
value function can be described by an appropriate
Hamilton--Jacobi--Bellman (HJB) equation. As time tends to infinity, a
heuristic separation of time and space variables in the HJB equation
yields an \textit{ergodic Bellman equation}. The optimal growth rate and
an optimal trading strategy are characterized by a specific solution of
this ergodic Bellman equation.

In contrast to (\ref{eq:gl354}), our robust problem (\ref{eq:gl532})
involves also the minimization among the class $\mathcal{Q}$, and this
would lead to a \textit{stochastic differential game} on an infinite time
horizon. Our main purpose, however, is to develop an alternative
approach: the main idea consists of combining the duality approach in~\cite{schiedwu} with methods from risk-sensitive control. Our main
results characterize the optimal growth rate
\[
\Lambda(\lambda):=\sup_{\xi}\mathop{\overline
{\lim}}_{T\uparrow\infty} \frac{1}{T}\ln\inf_{Q\in\mathcal
{Q}}E_Q[(X_T^{x_0,\xi})^\lambda],
\]
an optimal long-term investment strategy and an asymptotic worst-case
model $Q^*\in\mathcal{Q}$ for robust expected power utility in terms of
an appropriate ergodic Bellman equation.

Such asymptotic results on robust utility maximization are not only of
intrinsic interest but also relevant in connection to ``robust
large deviations'' criteria to optimal long term
investment. Suppose that the investor takes into account a class
$\mathcal{Q}$ of prior models and wants to maximize the worst-case
probability that the portfolio's growth rate $L_T^{x_0,\xi}:=\frac
{1}{T}\ln X_T^{x_0,\xi}$ exceeds some threshold $c\in\mathbb{R}$. In
the spirit of large deviations theory (see, e.g., \cite
{dembozeitouni}) the asymptotic problem then consists of
%
\begin{equation}
\label{eq:gl533}
\mbox{maximizing } \mathop{\overline{\lim}}_{T\uparrow\infty
} \frac{1}{T}\ln\inf_{Q\in\mathcal{Q}} Q[L_T^{x_0,\xi}\geq
c] \mbox{ among all } \xi.
\end{equation}
The solution can be derived by a duality approach similar to the G\"
{a}rtner--Ellis theorem, but here the dual problem involves the optimal
growth ra\-tes~$\Lambda(\lambda)$, $\lambda\in(0,1)$, of robust expected
power utility.

The paper is organized as follows: the setup is introduced in Section
\ref{sec:utilitymax}. Section~\ref{sec:outline} contains a heuristic
derivation of our main results that are verified in Section \ref
{sec:verification}. In Section \ref{sec:existenceergodicBellman} we
discuss the existence of a solution to our ergodic Bellman equation.
Explicit case studies are given in Section \ref{sec:explicitresults}.
In Section~\ref{sec:outperformance} we describe the duality approach to
the robust outperformance criterion~(\ref{eq:gl533}).

\section{The model and problem formulation}
\label{sec:utilitymax}
Let $(\Omega,\mathcal{F},(\mathcal{F}_t)_{t\geq0},Q_0)$ be the
canonical path space of a two-dimensional Wiener process
$W=(W^1_t,W^2_t)_{t\geq0}$.
We shall consider a long-term horizon investment model with one locally
riskless asset $S^0$ and one risky asset $S^1$. The performance of the
market is determined by an external ``economic factor''
$Y$, driven by the Wiener process $W$. The spectrum of possible
factors includes dividend yields, short-term interest rates,
price-earning ratios, yields on various bonds, the rate of inflation,
etc$\ldots.$ Both the price processes $S^0$, $S^1$ and the factor process~$Y$
will be subject to model ambiguity. This will be described by a class
$\mathcal{Q}$ of probabilistic models, viewed as perturbations of the
following reference model $Q_0$. Under $Q_0$ the dynamics of the
locally riskless asset is given by
\[
dS^0_t=S^0_tr(Y_t)\dd t,\qquad S^0_0=1,
\]
and the price process of the risky asset is governed by the SDE
%
\begin{equation}
\label{eq:gl333}
dS^1_t=S^1_t\bigl(m(Y_t)\dd t+\sigma\dd W^1_t\bigr).
\end{equation}
Thus the market price of risk is defined by
%
\begin{equation}
\label{eq:gl334}
\theta(y):=\frac{m(y)-r(y)}{\sigma}.
\end{equation}
The factor process evolves according to
%
\begin{equation}
\label{eq:gl274}
dY_t=g(Y_t)\dd t+\rho\dd W_t=g(Y_t)\dd t+\rho_1\dd W^1_t+\rho_2\dd W^2_t.
\end{equation}
We suppose that the economic factor can be observed but cannot be
traded directly. Therefore the market model is typically incomplete.
This class of market models is widely used in mathematical finance and
economics (see, e.g., \cite{fouque,fleminghernandez,castanedahernandez2}
and the references therein). Typically the
diffusion $Y$ is also assumed to be mean reverting and ergodic with
some invariant distribution~$\mu$. A special example is the
Ornstein--Uhlenbeck (OU) process with dynamics
%
\begin{equation}
\label{eq:gl420}
dY_t=\eta_0(\overline{y}-Y_t)\dd t+\sigma\dd W^1_t, \qquad\eta_0>0,\sigma
\not=0,
\end{equation}
and invariant distribution $\mu=\mathrm{N}(\overline{y},\frac{\sigma
^2}{2\eta_0})$.

We shall use the following general assumptions on the coefficients of
the diffusions, summarized as
\begin{assumption}
\label{th:assump1}
The functions $g$, $m$ admit derivatives $g_y,m_y\in C_b^1(\mathbb
{R})$, and $r$ belongs to $C_b^2(\mathbb{R})$, where $C_b^k(\mathbb
{R})$ denotes the class of all bounded functions with bounded
derivatives up to order $k$. Moreover, we assume that $\sigma$ and
$\Vert\rho\Vert$ are positive and that the short-rate function $r$ is
bounded below by some constant $a_1>0$.
\end{assumption}

Here we use \mbox{$\Vert\cdot\Vert$} to indicate the Euclidian norm in $\mathbb
{R}^2$, and in the sequel $(\cdot,\cdot)$ will denote the corresponding
inner product. In particular, our assumptions ensure\vadjust{\goodbreak} that the functions
$g$ and $\theta$ satisfy the linear growth conditions\looseness=1
\[
\vert g(y)\vert\leq a_2(1+\vert y\vert) \quad\mbox{and}\quad \vert\theta(y)\vert
\leq a_3\vert y\vert+a_4 \qquad\mbox{for } a_2,a_3,a_4>0.
\]\looseness=0
Note also that Assumption \ref{th:assump1} is consistent with
linear drift functions $g$ and~$m$. In this paper, such a choice of the
reference model will be particularly useful to obtain explicit
solutions (cf. Section \ref{sec:explicitresults}).

In reality, however, the ``true'' price dynamics
are not really known exactly. Here we focus on model uncertainty with
respect to the drift terms appearing in (\ref{eq:gl333}) and (\ref
{eq:gl274}). More precisely, we consider the parameterized class of
possible probabilistic models
\[
\mathcal{Q}:=\{Q^{\eta}|\eta=(\eta_t)_{t\geq0}\in\mathcal{C}\}
\]
on $(\Omega,\mathcal{F})$, where $\mathcal{C}$ denotes the set of all
progressively measurable processes $\eta=(\eta_t)_{t\geq0}$ such that
$\eta_t=(\eta^{11}_t,\eta^{12}_t,\eta^{21}_t,\eta^{22}_t)$ belongs
$dt\otimes Q_0$-a.e. to some fixed compact and convex set $\Gamma
\subset\mathbb{R}^4$ which contains the origin. For $\eta\in\mathcal
{C}$ and any fixed horizon $T$, the restriction of $Q^{\eta}$ to the
$\sigma$-field $\mathcal{F}_T$ is given by the Radon--Nikod\'{y}m density,
%
\begin{equation}
\label{eq:gl381}
D_T^{\eta}:=\frac{dQ^{\eta}}{dQ_0}\bigg|_{\mathcal{F}_T}:=\mathcal
{E}\biggl(\int_0^{\cdot}\eta^{1\cdot}_{t}Y_t+\eta^{2\cdot}_t\dd W_t\biggr)_T
\end{equation}
with respect to the reference measure $Q_0$. Here $\mathcal{E}(\cdot)$
denotes the It\^{o} exponential.
To see that $D_T^{\eta}$ is indeed the density of a probability measure
on~$(\Omega,\mathcal{F}_T)$, we can argue as follows: by
Assumption \ref{th:assump1} the diffusion process~$Y$ satisfies the
regularity conditions required in Lemma \ref{th:liptserersatz},
and so there exists some $\delta>0$ such that $\sup_{0\leq t\leq
T}E_{Q_0}[\exp(\delta Y_t^2)]<\infty$. The compactness of $\Gamma$ thus
ensures that
%
\begin{equation}
\label{eq:gl480}
\sup_{0\leq t\leq T}E_{Q_0}[\exp(\varepsilon\Vert\eta^{1\cdot
}_{t}Y_t+\eta
^{2\cdot}_t\Vert^2)]<\infty
\end{equation}
as soon as $\varepsilon>0$ is chosen sufficiently small. According to
\cite{liptsershiryaev}, Example 3, Section 6.2, this yields
$E_{Q_0}[D_T^{\eta}]=1$ as desired.

In view of (\ref{eq:gl381}) we have $Q_0=Q^{0}\in\mathcal{Q}$, and it
follows as in \cite{hernandezschied}, Lemma~3.1, that $\mathcal
{Q}$ is a convex set of locally equivalent measures on $(\Omega,\mathcal
{F})$. By Girsanov's theorem,
\[
W_t^{\eta}:=\biggl(W^1_t-\int_0^t\eta^{11}_sY_s+\eta^{21}_s\dd s,W^2_t-\int
_0^t\eta^{12}_sY_s+\eta^{22}_s\dd s\biggr),\qquad t\geq0,
\]
is a two-dimensional Wiener process under the measure $Q^{\eta}$, and
the dynamics of $S^1$, $Y$ under $Q^{\eta}$ take the form
%
\begin{subequation}
\label{eq:gl305}
\begin{eqnarray}\label{eq:gl267}
dY_t&=&[g(Y_t)+(\rho,\eta^{1\cdot}_tY_t+\eta^{2\cdot}_t)]\dd t+\rho\dd
W_t^{\eta},\\
dS^1_t&=&S^1_t\bigl([m(Y_t)+\sigma(\eta
^{11}_tY_t+\eta^{21}_t)]\dd t+\sigma\dd W_t^{1,\eta}\bigr).
\end{eqnarray}
\end{subequation}
Roughly speaking each element of $\mathcal{Q}$ corresponds to an \textit
{affine perturbation} of the drifts in our reference model $Q_0$. In
particular, our ``robust'' market model includes
the following special cases (see Section \ref{sec:explicitresults}):
\begin{example}[(Black--Scholes model with uncertain drift)]
\[
r(y)\equiv r,\qquad m(y)\equiv m,\qquad \Gamma=\{
(0,0) \}\times[a,b]\times\{ 0 \}.
\]
\end{example}
\begin{example}[(Geometric OU model with uncertain mean reversion)]
The factor process $Y$ is an OU process under $Q_0$ with rate of mean
reversion $\eta_0>0$, mean reversion level $\overline{y}=0$ and
volatility $\sigma>0$ [cf. (\ref{eq:gl420})]. We also assume $S^0_t=\exp
(rt)$, $r>0$ and $S_t^1:=\exp(Y_t+\alpha t)$, $\alpha\in\mathbb{R}$. By
It\^{o}'s formula this corresponds to
\[
g(y)=-\eta_0y,\qquad \rho_1=\sigma,\qquad \rho_2= 0,\qquad m(y)=-\eta_0y+\tfrac
{1}{2}\sigma^2+\alpha.
\]
Moreover, we take the set $\Gamma:=[\frac{\eta_0-b}{\sigma},\frac{\eta
_0-a}{\sigma}]\times\{(0,0,0)\}$ for $0<a\leq b<\infty$. For any
$Q^{\eta}\in\mathcal{Q}$ the process $Y$ thus follows under $Q^{\eta}\in
\mathcal{Q}$ OU-type dynamics with mean reversion process $\eta_0-\sigma
\eta_t^{11}$, taking values in $[a,b]$.
\end{example}

Let us now formulate our main problem. We consider an investor with
initial capital $x_0>0$ who aims at optimizing his portfolio in the
long run. A \textit{trading strategy} will be a predictable stochastic
process $\xi=(\xi^0,\xi^1)$ whose components $\xi^0$ and $\xi^1$
describe the successive amounts invested into the bond and into the
risky asset. The value of such a portfolio at time $t$ is given by
$X^\xi_t=\xi^0_t S^0_t+\xi^1_tS^1_t$. We also assume that $\xi^1$ is
$S^1$-integrable. Such a~trading strategy $\xi$ is said to be \textit
{self-financing} for the given initial capital~$x_0$ if its wealth
process $X^\xi=(X^\xi_t)_{t\geq0}$ takes the form
%
\begin{equation}
\label{eq:gl400}
X^\xi_t=x_0+\int_0^t\xi^0_u\dd S^0_u+\int_0^t\xi^1_u\dd S^1_u.
\end{equation}
Here the (stochastic) integrals can be interpreted as cumulative gains
or losses, that is, any change in the portfolio value equals the profit
or loss due to changes in the asset prices. For notational convenience
we omit the explicit dependence of $X^{\xi}$ on the initial capital
$x_0$, since it will be irrelevant for our purpose of long-term investment.
\begin{definition}
\label{th:definitionstrategyadmissible}
A self-financing trading strategy $\xi$ is called \textit
{$T$-admissi\-ble} if $X_t^\xi\geq0$ for all $t\in[0,T]$. A strategy $\xi
$ will be called admissible if it is $T$-admissible for any time
horizon $T>0$. We denote by $\mathcal{A}_T$ the class of all
$T$-admissible strategies and by $\mathcal{A}$ the class of all
admissible strategies.
\end{definition}

Clearly, a self-financing trading strategy $\xi$ can also be described
by the fractions
\[
\pi_t:=\frac{\xi^1_tS^1_t}{X_t^\xi},\qquad t\geq0,\vadjust{\goodbreak}
\]
of the current wealth which should be invested into the risky asset.
Throughout this paper we identify a strategy $\xi$ with the fractions
$\pi=(\pi_t)_{t\geq0}$. In terms of $\pi$ the wealth process defined
in (\ref{eq:gl400}) takes the form
\[
X^\pi_t=x_0+\int_0^t\frac{X_u^\pi(1-\pi_u)}{S^0_u}\dd S^0_u+\int
_0^t\frac{X^\pi_u\pi_u}{S^1_u}\dd S^1_u;
\]
that is, the investor's wealth $X^\pi$ evolves according to the SDE
%
\begin{eqnarray}
\label{eq:gl89}
dX^{\pi}_t&=&X^{\pi}_t\biggl((1-\pi_t)\frac{dS^0_t}{S^0_t}+\pi_t\frac
{dS^1_t}{S^1_t}\biggr)\nonumber\\[-8pt]\\[-8pt]
&=&X^{\pi}_t\bigl(r(Y_t)\dd t+\pi_t\sigma\bigl[\bigl(\theta(Y_t)+\eta^{11}_tY_t+\eta
^{21}_t\bigr)\dd t+dW^{1,\eta}_t\bigr]\bigr)\nonumber
\end{eqnarray}
with initial condition $X_0^{\pi}=x_0$.

In order to specify optimality, we assume that the investor's
preferences in the face of model ambiguity are described by a power
utility function
\[
u(x)=\frac{1}{\lambda}x^\lambda\qquad\mbox{with risk aversion
parameter } \lambda\in(0,1),
\]
and the set of prior probabilistic models $\mathcal{Q}$ (cf. page
\pageref{eq:gl150}). For a finite maturity~$T$, his robust portfolio
selection problem then consists of
%
\begin{equation}
\label{eq:gl295}
\mbox{maximizing } \inf_{Q^\eta\in\mathcal{Q}}E_{Q^\eta}[u(X_T^{\pi
})] \mbox{ among all } \pi\in\mathcal{A}_T.
\end{equation}
In a general semimartingale setting, this problem is well understood
from a theoretical point of view, in particular due to the articles
\cite{quenez,schiedwu,foellmergundel}. For robust
market models of the diffusion type described above and for power
utility, problem (\ref{eq:gl295}) has been discussed recently by Schied
\cite{schied1}. Applying dynamic programming methods to the dual
problem, he determines the maximal robust expected utility and a
worst-case model in terms of a~\textit{Hamilton--Jacobi--Bellman
equation}. Here we do not limit the analysis to a~fixed maturity.
Instead the objective of our investor consists of maximizing the
long-term growth of robust expected power utility. A~priori estimates,
as established in Lemma \ref{th:apriori}, suggest that the
maximal values
%
\begin{equation}
\label{eq:gl411}\qquad
U_T^Q(x_0):=\sup_{\pi\in\mathcal{A}_T}E_Q[u(X_T^{\pi})],\qquad
U_T(x_0):=\sup_{\pi\in\mathcal{A}_T}\inf_{Q^{\eta}\in\mathcal
{Q}}E_{Q^{\eta}}[u(X_T^{\pi})]
\end{equation}
for the classical utility maximization problem under $Q$ and for
its\break
robust extension will grow exponentially as $T\uparrow\infty$. Thus it
is natural to try to
%
\begin{equation}
\label{eq:gl27}
\mbox{maximize } \mathop{\overline{\lim}}_{T\uparrow
\infty} \frac{1}{T}\ln\inf_{Q^{\eta}\in\mathcal{Q}}E_{Q^{\eta
}}[(X_T^{\pi})^\lambda]
\mbox{ among all } \pi\in
\mathcal{A}.
\end{equation}
Our goal is to identify the optimal growth rate,
%
\begin{equation}
\label{eq:gl405}
\Lambda(\lambda):=\sup_{\pi\in\mathcal{A}}\mathop{\overline
{\lim}}_{T\uparrow\infty} \frac{1}{T}\ln\inf_{Q^{\eta}\in\mathcal
{Q}}E_{Q^{\eta}}[(X_T^{\pi})^\lambda], \qquad\lambda\in(0,1),
\end{equation}
an optimal long term investment strategy $\pi^*$ and an asymptotic
worst-case model $Q^{\eta^*}\in\mathcal{Q}$. Heuristically this means
that, as $T\uparrow\infty$,
%
\begin{eqnarray}\label{eq:gl337}
U_T(x_0)&\approx&\frac{1}{\lambda}x_0^\lambda e^{\Lambda(\lambda)
T}\\
\label{eq:gl340}
&\approx&
\inf_{Q^{\eta}\in\mathcal{Q}}E_{Q^{\eta}}[ u(X_T^{\pi^*})]\\
\label{eq:gl341}
&\approx& U_T^{Q^{\eta^*}}(x_0)=\sup_{\pi\in
\mathcal{A}_T}E_{Q^{\eta^*}}[u(X_T^{\pi})]\\
\label{eq:gl342}
&\approx&
E_{Q^{\eta^*}}[u(X_T^{\pi^*})].
\end{eqnarray}
Here (\ref{eq:gl340}) corresponds to asymptotic optimality of the
trading strategy~$\pi^*$, (\ref{eq:gl341}) to the property of $Q^{\eta
^*}$ of being the asymptotic worst-case model, and~(\ref{eq:gl342})
identifies $\pi^*$ also as the asymptotically optimal strategy for the
model~$Q^{\eta^*}$. In particular,~$Q^{\eta^*}$ and~$\pi^*$ can be
viewed as a saddle point for the problem of asymptotic robust utility
maximization with control parameters $\eta\in\mathcal{C}$ and $\pi\in
\mathcal{A}$. Moreover, (\ref{eq:gl340})~suggests that an optimal
strategy $\pi^*$ of the asymptotic criterion (\ref{eq:gl27}) should
provide a good approximation of an optimal investment process $\pi
^{*,T}$ for the robust power utility maximization problem with a large
but finite time horizon $T$.
\begin{remark}
\label{th:negpow}
The asymptotic approach to robust utility maximization can be extended
to the following cases (see \cite{knispel}, Chapter 4):
\begin{itemize}
\item
For power utility $u(x)=\frac{1}{\lambda}x^\lambda$ with parameter
$\lambda<0$ the distance between
\[
U_T(x_0)=\frac{1}{\lambda}\inf_{\pi\in\mathcal{A}_T}\sup_{Q^{\eta}\in
\mathcal{Q}} E_{Q^{\eta}}[(X_T^{\pi})^\lambda]
\]
and its upper bound $0$ will typically decrease exponentially as
$T\uparrow\infty$. This suggests that we should compute the optimal
growth rate,
\[
\Lambda(\lambda):=\inf_{\pi\in\mathcal{A}}\mathop{\underline{\lim}}_{T\uparrow\infty} \frac{1}{T} \ln\sup_{Q^{\eta}\in
\mathcal{Q}} E_{Q^{\eta}}[(X_T^{\pi})^\lambda].
\]
\item
For logarithmic utility $u(x)=\ln(x)$ the growth of robust expected
utility will be linear. Thus we want to
\[
\mbox{maximize }
\mathop{\overline{\lim}}_{T\uparrow\infty} \frac{1}{T}\inf
_{Q^{\eta}\in\mathcal{Q}}E_{Q^{\eta}}[\ln(X_T^{\pi})] \mbox{ among
all } \pi\in\mathcal{A}.
\]
\end{itemize}
\end{remark}
%
\section{Heuristic outline of the dynamic programming approach}
\label{sec:outline}
We start with a heuristic derivation of our main results. They provide
a characterization of the optimal growth rate $\Lambda(\lambda)$, of an
asymptotic worst-case model~$Q^{\eta^*}$, and of an optimal long-term
investment strategy $\pi^*$ in terms of an \textit{ergodic Bellman
equation} (EBE). Our method combines the \textit{duality approach} to
robust utility maximization with \textit{dynamic programming methods} for
a varying time horizon. As a byproduct of the duality approach, we also
show that~$U_T(x_0)$ grows exponentially at rate $\Lambda(\lambda)$ as
$T\uparrow\infty$. A more direct, but not more tractable approach to
the saddle-point problem (\ref{eq:gl27}) via stochastic differential
games will be discussed in Remark \ref{th:differentialgame}.

First, we set up the duality approach based on the results of Schied and
Wu \cite{schiedwu} for a utility function $u$ on the positive halfline.
This will allow us to transform the primal saddle-point problem (\ref
{eq:gl27}) to a simpler minimization problem on the dual side. The dual
value function at time $T$ is defined by
%
\begin{equation}
\label{eq:gl504}
V_T(y):=\inf_{Q\in\mathcal{Q}}\inf_{Y\in\mathcal
{Y}^Q}E_Q[v(yY_T/S^0_T)],\qquad y>0,
\end{equation}
where $v(y):=\sup_{x>0}\{u(x)-xy\}$, $y>0$, is the convex conjugate
function of~$u$. This definition also involves the class of supermartingales
\[
\mathcal{Y}_T^{Q}:=\{Y\geq0|Y_0=1 \mbox{ and }\forall\pi\in\mathcal
{A}_T\dvtx(Y_tX_t^{\pi}/S^0_t)_{t\leq T}\mbox
{ is a }Q\mbox{-supermartingale}\}
\]
as introduced by Kramkov and Schachermayer \cite{kramkovschachermayer}.
Note that $\mathcal{Y}_T^Q$ contains the density processes (taken with
respect to $Q$ and the num\'{e}raire $S^0$) of the class $\mathcal
{P}_T$ of all equivalent local martingale measures on $(\Omega,\mathcal
{F}_T)$. For power utility we have $v(y)=-\beta^{-1}y^\beta$, $\beta
:=\frac{\lambda}{\lambda-1}$, and this yields the scaling property
$V_T(y)=y^\beta V_T(1)$. Due to \cite{schiedwu}, Theorem 2.2, the
primal value function (\ref{eq:gl411}) can then be obtained as
%
\begin{equation}
\label{eq:gl506}
U_T(x_0)=\inf_{y>0}\{V_T(y)+ x_0y\}=\frac{1}{\lambda}x_0^\lambda
(-\beta V_T(1))^{1-\lambda}.
\end{equation}
Since power utility has asymptotic elasticity $\overline
{\lim}_{x\uparrow\infty}\frac{xu'(x)}{u(x)}<1$, it follows from \cite
{schiedwu}, Theorem 2.5, also that
%
\begin{equation}
\label{eq:gl505}
V_T(1)=\inf_{P\in\mathcal{P}_T}\inf_{Q\in\mathcal{Q}}E_Q\biggl[v\biggl(\frac
{dP}{dQ}\bigg|_{\mathcal{F}_T}\Big/S^0_T\biggr)\biggr].
\end{equation}
We now parameterize the sets $\mathcal{Y}_T^Q$ and $\mathcal{P}_T$.
Since $Z_t:=dP/dQ_0|_{\mathcal{F}_t}$, $t\leq T$, is a positive
$Q_0$-martingale for any $P\in\mathcal{P}_T$, the martingale
representation theorem yields the existence of an $\mathbb{R}^2$-valued
progressively measurable process $\phi=(\phi^1,\phi^2)$ with $\int
_0^T\Vert\phi_s\Vert^2\dd s<\infty$ $Q_0$-a.s. such that $Z_t=\mathcal
{E}(\int_0^{\cdot}\phi_s\dd W_s)_t$. By Girsanov's theorem, the
discounted wealth process $X^\pi/S^0$ is a local martingale under $P$
if and only if $\phi^1_s=-\theta(Y_s)$ $ds\otimes Q_0$-a.e. Thus the
$Q_0$-density process of an martingale measure $P\in\mathcal{P}_T$
necessarily takes the form
%
\begin{equation}
\label{eq:gl384}
Z^{\nu}_t:=\mathcal{E}\biggl(-\int_0^{\cdot}\theta(Y_s)\dd W_s^1-\int_0^{\cdot
}\nu_s\dd W_s^2\biggr)_t
\end{equation}
for some progressively measurable process $\nu$ such that $\int_0^T\nu
^2_s\dd s<\infty$ $Q_0$-a.s. Conversely, $Z_T^{\nu}$ corresponds to
the $Q_0$-density of an equivalent local martingale measure on $(\Omega
,\mathcal{F}_T)$ as soon as the martingale condition $E_{Q_0}[Z_T^{\nu
}]=1$ holds. This can be verified if, for instance, the process $\nu$ is
assumed to be bounded. Thus our market model admits a variety of
equivalent local martingale measures up to any finite horizon $T$; that
is, the restriction of our model to a finite horizon is \textit
{arbitrage-free} but \textit{incomplete}.

More generally, we will denote by $\mathcal{M}$ the set of all
progressively measurable processes $\nu=(\nu_t)_{t\geq0}$ such that
$\int_0^T\nu^2_t\dd t<\infty$ $Q_0$-a.s. for all $T>0$. Via~(\ref
{eq:gl384}) every $\nu\in\mathcal{M}$ gives rise to a positive
$Q_0$-supermartingale $Z^{\nu}$. Using It\^{o}'s formula one easily
shows that $(D^{\eta})^{-1}Z^\nu X^\pi/S^0$ is a positive local
martingale under $Q^{\eta}$ for any $\nu\in\mathcal{M}$ and $\pi\in
\mathcal{A}_T$, and hence a $Q^\eta$-supermartingale. Thus
\[
\biggl\{\biggl(\frac{dP}{dQ^{\eta}}\bigg|_{\mathcal{F}_t}\biggr)_{t\leq T}\Big|P\in\mathcal
{P}_T\biggr\}\subset\{((D_t^{\eta})^{-1}Z_t^{\nu})_{t\leq T}|\nu\in\mathcal
{M}\}\subset\mathcal{Y}_T^{Q^{\eta}}.
\]
In view of (\ref{eq:gl504}), (\ref{eq:gl505}) and (\ref{eq:gl506}) this
inclusion and a change of measure yield
%
\begin{equation}
\label{eq:gl47}
U_T(x_0)=\frac{1}{\lambda}x_0^\lambda\Bigl(\inf_{\nu\in\mathcal{M}}\inf
_{\eta\in\mathcal{C}}E_{Q_0}\bigl[(Z_T^{\nu}(S^0_T)^{-1})^{{\lambda
}/({\lambda-1})} (D_T^{\eta})^{{1}/({1-\lambda})}\bigr]\Bigr)^{1-\lambda}.
\end{equation}
In a second step, we derive an ergodic Bellman equation by applying
dynamic programming methods to the dual minimization problem. Since
$Z_T^{\nu}$, $D_T^{\eta}$ and the bond price $S^0_T$ depend on the
factor process $Y$, the expectation at the right-hand side of (\ref
{eq:gl47}) is a function of the initial state $Y_0=y$. For all
processes $\eta\in\mathcal{C}$ and $\nu\in\mathcal{M}$ we can thus define
%
\begin{equation}
\label{eq:gl37}
V(\eta,\nu,y,T):=E_{Q_0}\bigl[(Z_T^{\nu}(S^0_T)^{-1})^{{\lambda
}/({\lambda-1})} (D_T^{\eta})^{{1}/({1-\lambda})}\bigr].
\end{equation}
Inserting the definitions of $Z_T^{\nu}$, $D_T^{\eta}$ and $S^0_T$ we
then obtain the decomposition
%
\begin{equation}
\label{eq:gl507}
V(\eta,\nu,y,T)=E_{Q_0}\bigl[\mathcal{E}^{\eta,\nu}_T e^{\int_0^Tl(\eta_t,\nu
_t,Y_t)\dd t}\bigr].
\end{equation}
Here the function $l\dvtx\Gamma\times\mathbb{R}\times\mathbb{R}\rightarrow
\mathbb{R}_+$ is defined by
%
\begin{eqnarray}
\label{eq:gl135} l(\eta,\nu,y)&:=&\frac{1}{2}\frac{\lambda}{(1-\lambda
)^2}\bigl[\bigl(\theta(y)+\eta^{11}y+\eta^{21}\bigr)^2+(\nu+\eta^{12}y+\eta
^{22})^2\bigr]\nonumber\\[-8pt]\\[-8pt]
&&{}+\frac{\lambda}{1-\lambda}r(y)\nonumber
\end{eqnarray}
and
\[
\mathcal{E}_T^{\eta,\nu}:=\mathcal{E}\biggl(\frac{1}{1-\lambda}\biggl(\int_0^{\cdot
}\lambda\theta(Y_t)+\eta^{11}_tY_t+\eta^{21}_t\dd W^1_t+\int_0^{\cdot
}\lambda\nu_t+\eta^{12}_tY_t+\eta^{22}_t\dd W^2_t\biggr)\biggr)_T.
\]
To simplify the expression for $V(\eta,\nu,y,T)$, we shall interpret
the It\^{o} exponential as the density of a probability measure $R^{\eta
,\nu}$ on $(\Omega,\mathcal{F}_T)$. This requires $E_{Q_0}[\mathcal
{E}_T^{\eta,\nu}]=1$ which is satisfied, for example, if $\int_0^T\nu
_t^2\dd t$ is bounded. For arbitrary $\nu\in\mathcal{M}$ we may have
$E_{Q_0}[\mathcal{E}_T^{\eta,\nu}]<1$, but here we argue heuristically,
and so we postpone this technical problem to the proof of Theorem
\ref{th:verification1}. In terms of the measure $R^{\eta,\nu}$ we can write
%
\begin{equation}
\label{eq:gl138}
V(\eta,\nu,y,T)=E_{R^{\eta,\nu}}\bigl[e^{\int_0^Tl(\eta_t,\nu_t,Y_t)\dd t}\bigr].
\end{equation}
Moreover, Girsanov's theorem yields that the factor process
$(Y_t)_{t\leq T}$ evolves under $R^{\eta,\nu}$ according to the SDE
%
\begin{equation}
\label{eq:gl44}
dY_t=h(\eta_t,\nu_t,Y_t)\dd t+\rho\dd W_t^{\eta,\nu},
\end{equation}
where $W^{\eta,\nu}$ is a Wiener process under $R^{\eta,\nu}$ and where
$h$ is defined by
%
\begin{eqnarray}
\label{eq:gl42}
h(\eta,\nu,y)&:=&g(y)+\frac{1}{1-\lambda}\rho_1\bigl(\lambda\theta(y)+\eta
^{11}y+\eta^{21}\bigr)\nonumber\\[-8pt]\\[-8pt]
&&{}+\frac{1}{1-\lambda}\rho_2(\lambda\nu+\eta
^{12}y+\eta^{22}).
\nonumber
\end{eqnarray}
Putting (\ref{eq:gl47}), (\ref{eq:gl37}) and (\ref{eq:gl138}) together,
we get
%
\begin{equation}
\label{eq:gl335}
U_T(x_0)=\frac{1}{\lambda}x_0^\lambda v(y,T)^{1-\lambda},
\end{equation}
where
\[
v(y,T):=\inf_{\nu\in\mathcal{M}}\inf_{\eta\in\mathcal{C}}E_{R^{\eta,\nu
}}\bigl[e^{\int_0^Tl(\eta_t,\nu_t,Y_t)\dd t}\bigr]
\]
denotes the value function of the finite horizon optimization problem
on the dual side of (\ref{eq:gl47}). Such an ``expected
exponential of integral cost criterion'' with a
dynamics of the form (\ref{eq:gl44}) is standard in stochastic control
theory (see, e.g., \cite{flemingsoner}, Remark IV.3.3). As a
result, $v$ can be described as the solution to the \textit
{Hamilton--Jacobi--Bellman} (HJB) equation,
%
\begin{equation}
\label{eq:gl23}
v_t=\frac{1}{2}\Vert\rho\Vert^2 v_{yy}+\inf_{\nu\in\mathbb{R}}\inf
_{\eta\in\Gamma}\{l(\eta,\nu,\cdot)v +h(\eta,\nu,\cdot)v_y\},\qquad
v(\cdot,0)\equiv1.
\end{equation}

The following lemma establishes a priori bounds for the exponential
growth of robust expected power utility, and this justifies the scaling
in (\ref{eq:gl27}).
\begin{lemma}
\label{th:apriori}
Suppose in addition to Assumption \ref{th:assump1} that one of
the following conditions is satisfied:
\begin{longlist}[(2)]
\item[(1)] The market price of risk function $\theta$ in (\ref
{eq:gl334}) is bounded.
\item[(2)] There exist constants $K,M_1,M_2>0$ such that
\[
-K y+M_1\leq g(y)+\frac{\lambda}{1-\lambda}\rho_1\theta(y)\leq-K
y+M_2,\qquad 2\frac{\lambda}{(1-\lambda)^2}\Vert\rho\Vert^2 a_3^2<K^2.
\]
\end{longlist}
Then there are constants $K_1,K_2>0$ such that for any initial capital $x_0>0$
%
\begin{equation}
\label{eq:gl181}
K_1\leq\mathop{\underline{\lim}}_{T\uparrow\infty} \frac
{1}{T}\ln U_T(x_0)\leq\mathop{\overline{\lim}}_{T\uparrow\infty
} \frac{1}{T}\ln U_T(x_0)\leq K_2.
\end{equation}
\end{lemma}
\begin{pf}
If at any time the whole capital is put into the money market account,
then the investor's utility at time $T$ is given by $\frac
{1}{\lambda}x_0^\lambda\exp(\lambda\int_0^T r(Y_t)\dd t)$ which, by
Assumption \ref{th:assump1}, is bounded from below by $\frac
{1}{\lambda}x_0^\lambda\exp(\lambda a_1 T)$. This implies the lower bound
\[
0<K_1:=\lambda a_1\leq\mathop{\underline{\lim}}_{T\uparrow\infty
} \frac{1}{T}\ln U_T(x_0).
\]
To obtain the upper bound, observe first that
\[
v(y,T)\leq V(0,0,y,T)\leq E_R\bigl[e^{({1/2})({\lambda}/{(1-\lambda
)^2})\int_0^T\theta^2(Y_t)\dd t}\bigr]e^{({\lambda}/({1-\lambda}))\Vert
r\Vert
_\infty T},
\]
where $R:=R^{0,0}$ is the probability measure defined by $\mathcal
{E}^{0,0}_T$. In view of (\ref{eq:gl335}) we thus get the estimate
%
\begin{eqnarray}
\label{eq:gl336}
&&\mathop{\overline{\lim}}_{T\uparrow\infty} \frac{1}{T}\ln
U_T(x_0)\nonumber\\[-8pt]\\[-8pt]
&&\qquad\leq(1-\lambda) \mathop{\overline{\lim}}_{T\uparrow\infty
} \frac{1}{T}\ln E_{R}\bigl[e^{({1/2})({\lambda}/{(1-\lambda
)^2})\int_0^T\theta^2(Y_t)\dd t}\bigr] +\lambda\Vert
r\Vert_\infty.\nonumber
\end{eqnarray}
In particular, the upper bound in (\ref{eq:gl181}) holds with
$K_2:=\frac{1}{2}\frac{\lambda}{1-\lambda}\Vert\theta\Vert^2_\infty
+\lambda\Vert r\Vert_\infty$ if the market price of risk function
$\theta$ is bounded.

Case (2) requires more effort. By (\ref{eq:gl336}) it is sufficient to
show that
%
\begin{equation}
\label{eq:gl187}
\mathop{\overline{\lim}}_{T\uparrow\infty} \frac{1}{T}\ln
E_{R}\biggl[\exp\biggl(\frac{1}{2}\frac{\lambda}{(1-\lambda)^2}\int_0^T\theta
^2(Y_t)\dd t\biggr)\biggr] <\infty.
\end{equation}
To this end, recall from (\ref{eq:gl44}) that the dynamics of $Y$ under
$R$ are given by
\[
dY_t=h(0,0,Y_t)\dd t+\rho\dd W_t^{0,0} \qquad\mbox{with }
h(0,0,y)=g(y)+\frac{\lambda}{1-\lambda}\rho_1\theta(y).
\]
Consider now the $R$-OU processes $dZ_{it}=[-K Z_{it}+M_i]\dd t+\rho\dd
W_t^{0,0}$, $Z_{i0}=y$, $i=1,2$. Then a comparison argument for the
solutions of SDEs ensures that
%
\begin{equation}
\label{eq:gl185}
R[Z_{1t}\leq Y_t\leq Z_{2t} \mbox{ for all } t\geq0]=1.
\end{equation}
Take now $\varepsilon>0$ satisfying $2\frac{\lambda}{(1-\lambda
)^2}\Vert
\rho\Vert^2 (a_3^2+\varepsilon)<K^2$. By Assumption \ref{th:assump1}
there exist constants $C_1$, $C_2$ depending on $\varepsilon$ such that
\[
\theta^2(y)\leq(a_3\vert y\vert+a_4)^2\leq\biggl(a_3^2+\frac{\varepsilon
}{2}\biggr)y^2+C_1\leq(a_3^2+\varepsilon)(y-M_i/K)^2+C_1+C_2
\]
for any $y\in\mathbb{R}$. Together with (\ref{eq:gl185}) and H\"
{o}lder's inequality (applied in line 3) this leads to
%
\begin{eqnarray}
\label{eq:gl186}
&&E_{R}\bigl[e^{({1/2})({\lambda}/{(1-\lambda)^2})\int_0^T\theta
^2(Y_t)\dd t}\bigr]\nonumber\\
&&\qquad\leq E_{R}\bigl[e^{({1/2})({\lambda}/{(1-\lambda)^2})
(a_3^2+{\varepsilon}/{2}) \int_0^T Y^2_t\dd t}\bigr]e^{({1/2})
({\lambda}/{(1-\lambda)^2})C_1T}\nonumber\\
&&\qquad\leq E_{R}\bigl[e^{({1/2})({\lambda}/{(1-\lambda)^2})(a_3^2+
{\varepsilon}/{2}) \int_0^T Z^2_{1t}+Z^2_{2t}\dd
t}\bigr]e^{C_3T}\nonumber\\[-8pt]\\[-8pt]
&&\qquad\leq \max_{i=1,2} E_{R}\bigl[e^{{\lambda}/{(1-\lambda
)^2}(a_3^2+{\varepsilon}/{2}) \int_0^T Z^2_{it}\dd
t}\bigr]e^{C_3T}\nonumber\\
&&\qquad\leq \max_{i=1,2} E_{R}\bigl[e^{{\lambda}/{(1-\lambda
)^2}(a_3^2+\varepsilon) \int_0^T(Z_{it}-M_i/K)^2\dd t}\bigr]e^{C_4T}\nonumber\\
&&\qquad= \max_{i=1,2} E_{R}\bigl[e^{{\lambda}/{(1-\lambda)^2}\Vert\rho
\Vert^2(a_3^2+\varepsilon) \int_0^T \widetilde{Z}_{it}^2\dd
t}\bigr]e^{C_4T}.\nonumber
\end{eqnarray}
Here we use the processes $\widetilde{Z}_i$, $i=1,2$, defined by
$\widetilde{Z}_{it}:=\Vert\rho\Vert^{-1}(Z_{it}-M_i/K)$. Note that
$\widetilde{Z}_i$ is an OU process with rate of mean reversion $K$,
equilibrium level $0$ and volatility $1$, since $B:=\int_0(\Vert\rho
\Vert)^{-1}\rho\dd W^{0,0}_t$ is a standard one-dimensional $R$-Brownian
motion, due to L\'{e}vy's characterization.
Applying Lemma 4.2 in \cite{florenslandaispham} [here with
$\lambda=0$, $\mu=\frac{\lambda}{(1-\lambda)^2}\Vert\rho\Vert^2
(a_3^2+\varepsilon)$ and $\theta_0=-K$] for the asymptotics of the Laplace
transform of the energy integral of a~normalized OU process, we obtain
\begin{eqnarray*}
&&\lim_{T\uparrow\infty} \frac{1}{T}\ln E_{R}\bigl[e^{
({\lambda}/{(1-\lambda)^2})\Vert\rho\Vert^2(a_3^2+\varepsilon) \int
_0^T\widetilde{Z}_{it}^2\dd t}\bigr]\\
&&\qquad = \frac
{1}{2}\Biggl(K - \sqrt{K^2 - 2\frac{\lambda}{(1-\lambda)^2}\Vert\rho\Vert
^2 (a_3^2 + \varepsilon)} \Biggr).
\end{eqnarray*}
In view of (\ref{eq:gl186}) we have thus shown (\ref{eq:gl187}). This
completes the proof.
\end{pf}
Combining the discussion of (\ref{eq:gl337}) with (\ref{eq:gl335}), it
is natural to expect that the optimal growth rate $\Lambda(\lambda)$ in
(\ref{eq:gl405}) satisfies
\[
\Lambda(\lambda)=\mathop{\overline{\lim}}_{T\uparrow\infty}
\frac{1}{T}\ln U_T(x_0)=\mathop{\overline{\lim}}_{T\uparrow\infty
} \frac{1}{T}\ln(v(y,T)^{1-\lambda}).
\]
As in Fleming and McEneaney \cite{flemingmceneaney} we now use a formal
separation of time and space variables and formulate the heuristic ansatz
%
\begin{equation}
\label{eq:gl249}
(1-\lambda)\ln v(y,T)=\ln U_T(x_0)\approx\Lambda(\lambda)T+\varphi(y).
\end{equation}
Here the function $\varphi\dvtx\mathbb{R}\rightarrow\mathbb{R}$
incorporates the influence of the initial state $Y_0=y$. Inserting this
ansatz into the HJB equation (\ref{eq:gl23}), we obtain a steady-state
dynamic programming equation for the pair ($\Lambda(\lambda),\varphi$)
%
\begin{eqnarray}
\label{eq:gl39}
\Lambda(\lambda)&=&\frac{1}{2}\Vert\rho\Vert^2\biggl[\varphi_{yy}+\frac
{1}{1-\lambda}\varphi_y^2\biggr]\nonumber\\[-8pt]\\[-8pt]
&&{}+\inf_{\nu\in\mathbb{R}}\inf_{\eta\in\Gamma}\{
(1-\lambda)l(\eta,\nu,\cdot)+\varphi_yh(\eta,\nu,\cdot)\}.\nonumber
\end{eqnarray}
An equation of this type is called an \textit{ergodic Bellman equation}
(EBE) (see, e.g., \cite{bensoussanfrehse,kaisesheu,nagai}
and the references\vadjust{\goodbreak} therein). For fixed $\eta\in\Gamma$ the
minimi\-zer~$\nu^*(\eta,y)$ among all $\nu\in\mathbb{R}$ can be computed
explicitly as
%
\begin{equation}
\label{eq:gl390}
\nu^*(\eta,y)=-\eta^{12}y-\eta^{22}-\rho_2\varphi_y(y).
\end{equation}
Thus the EBE (\ref{eq:gl39}) can be rewritten in condensed form that
involves only an infimum among the set $\Gamma$. Let us now assume that
our EBE (\ref{eq:gl39}) admits a~solution $\Lambda(\lambda)\in\mathbb
{R}_+$, $\varphi\in C^2(\mathbb{R})$. In addition, assume that $\eta
^*(y)$ is a minimizer in (\ref{eq:gl39}), and let $Q^{\eta^*}\in\mathcal
{Q}$ be the probabilistic model corresponding to the feedback control
$\eta_t^*=\eta^*(Y_t)$. We are now going to give a heuristic argument
to identify a candidate for the optimal long-run investment process~$\pi
^*$. To this end, we suppose that the measure $Q^{\eta^*}$ is a
worst-case model in the asymptotic sense that
%
\begin{equation}
\label{eq:gl51}
\Lambda(\lambda)=\mathop{\overline{\lim}}_{T\uparrow\infty}
\frac{1}{T}\ln U_T(x_0)=\mathop{\overline{\lim}}_{T\uparrow\infty
} \frac{1}{T}\ln\sup_{\pi\in\mathcal{A}_T}E_{Q^{\eta^*}}[(X^{\pi
}_T)^\lambda].
\end{equation}
Later on we will show that this assumption is indeed justified. We are
now going to introduce a change of measure which will allow us to
interpret the finite time maximization problem at the right-hand side
of (\ref{eq:gl51}) as an exponential of integral criterion. For this
purpose, note that an optimal wealth process should stay positive, and
this suggests that we should focus on those strategies $\pi\in\mathcal
{A}$, where the unique strong solution to (\ref{eq:gl89}) takes the form
\[
X_t^{\pi}=x_0e^{\int_0^t\pi_u\sigma\dd W_u^{1,\eta^*}+\int
_0^tr(Y_u)+\sigma\pi_u(\theta(Y_u)+\eta_u^{11,*}Y_u+\eta
^{21,*}_u)-({1/2})\sigma^2\pi_u^2\dd u}.
\]
In this case the expectation at the right-hand side of (\ref{eq:gl51})
can be rewritten~as
\[
E_{Q^{\eta^*}}[(X_T^{\pi})^\lambda]=x_0^\lambda E_{R^{\pi,\eta
^*}}\bigl[e^{\int_0^T \widetilde{l}(\pi_t,\eta^*(Y_t),Y_t)\dd t}\bigr].
\]
Here we use the notation
%
\begin{equation}
\label{eq:gl108}
\widetilde{l}(\pi,\eta,y):=\tfrac{1}{2}\lambda(\lambda-1)\sigma^2\pi
^2+\lambda\sigma[\theta(y)+\eta^{11}y+\eta^{21}]\pi+\lambda r(y),
\end{equation}
and $R^{\pi,\eta}$ denotes the probability measure on $(\Omega,\mathcal
{F}_T)$ defined by
%
\begin{equation}
\label{eq:measuredifferential}
\frac{dR^{\pi,\eta}}{dQ^{\eta}}\bigg|_{\mathcal{F}_T}:=\mathcal{E}\biggl(\int
_0^{\cdot}\lambda\pi_t\sigma\dd W_t^{1,\eta}\biggr)_T.
\end{equation}
By Girsanov's theorem, the dynamics of $(Y_t)_{t\leq T}$ under
$R^{\pi,\eta}$ are described~by
%
\begin{equation}
\label{eq:gl338}
dY_t=\widetilde{h}(\pi_t,\eta_t,Y_t)\dd t+\rho\dd W_t^{\pi,\eta}
\end{equation}
in terms of the function $\widetilde{h}$ defined by
%
\begin{equation}
\label{eq:gl98}
\widetilde{h}(\pi,\eta,y):=g(y)+(\rho,\eta^{1\cdot}y+\eta^{2\cdot
})+\lambda\rho_1\sigma\pi
\end{equation}
and the one-dimensional Wiener process $W^{\pi,\eta}$. We have thus
shown that the finite horizon maximization problem appearing in the
right-hand side of~(\ref{eq:gl51}) can be viewed as a finite horizon
control problem with value function
\[
\widetilde{v}(y,T):=\sup_{\pi\in\mathcal{A}}E_{Q^{\eta^*}}[(X_T^{\pi
})^\lambda] =x_0^\lambda\sup_{\pi\in\mathcal{A}}E_{R^{\pi,\eta
^*}}\bigl[e^{\int_0^T \widetilde{l}(\pi_t,\eta^*(Y_t),Y_t)\dd t}\bigr]
\]
and with dynamics (\ref{eq:gl338}). In analogy to (\ref{eq:gl23}), we
expect that $\widetilde{v}$ is the solution to the HJB equation
%
\begin{equation}
\label{eq:gl339}
\widetilde{v}_t=\frac{1}{2}\Vert\rho\Vert^2 \widetilde{v}_{yy}+\sup
_{\pi\in\mathbb{R}}\{\widetilde{l}(\pi,\eta^*,\cdot)\widetilde{v}
+\widetilde{h}(\pi,\eta^*,\cdot)\widetilde{v}_y\}, \qquad\widetilde
{v}(\cdot,0)\equiv1.
\end{equation}
Our ansatz (\ref{eq:gl249}) combined with (\ref{eq:gl51}) for the
worst-case measure $Q^{\eta^*}$ now suggests the heuristic separation
of variables $\ln\widetilde{v}(y,T)\approx\Lambda(\lambda)T+\varphi
(y)$. Inserting this asymptotic identity into (\ref{eq:gl339}), we
finally obtain an alternative version of the EBE,
%
\begin{equation}
\label{eq:gl28}
\Lambda(\lambda)=\frac{1}{2}\Vert\rho\Vert^2[\varphi_{yy}+\varphi
^2_y]+\sup_{\pi\in\mathbb{R}}\{\widetilde{l}(\pi,\eta^*,\cdot)+\varphi
_y\widetilde{h}(\pi,\eta^*,\cdot)\}.
\end{equation}
Note that the role played by the controls $\eta$ and $\nu$ in (\ref
{eq:gl39}) is now taken over by the ``trading strategies''
$\pi$. We expect that the maximizing function
%
\begin{equation}
\label{eq:gl24}
\pi^*(y)=\frac{1}{1-\lambda}\frac{1}{\sigma}\bigl(\rho_1\varphi_y(y)+\theta
(y)+\eta^{11,*}(y)y+\eta^{12,*}(y)\bigr)
\end{equation}
in (\ref{eq:gl28}) provides an optimal feedback control $\pi^*_t=\pi
^*(Y_t)$, $t\geq0$, for the asymptotic maximization of power utility
with respect to the specific mo\-del~$Q^{\eta^*}$ and at the same time
for the original robust problem (\ref{eq:gl27}).

\section{Verification theorems}
\label{sec:verification}
In this section we verify our heuristic results. For this purpose, we
first return to the heuristic change of measure in (\ref{eq:gl138})
which is crucial to translate the dual problem (\ref{eq:gl47}) into a
standard\vspace*{1pt} ``exponential of integral criterion.'' From the
technical point of view this requires the condition $E_{Q_0}[\mathcal
{E}_T^{\eta,\nu}]=1$ that can be violated if the supermartingale $Z^{\nu
}$ is not a true $Q_0$-martingale. This fact will create some technical
difficulties. To overcome this obstacle, we shall employ a localization
argument.\looseness=1

\begin{lemma}
\label{th:localization1}
Let $\eta\in\mathcal{C}$ and $\nu\in\mathcal{M}$ be arbitrary controls,
and suppose that $(\tau_n)_{n\in\mathbb{N}}$ is a localizing sequence
of stopping times for the local $Q_0$-martingale~$Z^{\nu}$. Then $V(\eta
,\nu,y,T\wedge\tau_n)\nearrow V(\eta,\nu,y,T)$ as $n\uparrow\infty$,
and the integrands in (\ref{eq:gl37}) even converge in $L^1(Q_0)$ if
$V(\eta,\nu,y,T)<\infty$.
\end{lemma}

\begin{pf}
The proof is given in \cite{schied1}, Lemma 3.2. The main idea
consists of applying the concept of \textit{extended martingale measures}
introduced in \cite{foellmergundel}.
\end{pf}

In a second step we are going to show that the value $\widetilde{\Lambda
}(\lambda)$ given by a~specific solution to the EBE (\ref{eq:gl39}) is
actually the exponential growth rate of the maximal robust power
utility $U_T(x_0)$. For this purpose, we need
\begin{assumption}
\label{th:assump3}
Suppose that $\widetilde{\Lambda}(\lambda)\in\mathbb{R}_+$, $\varphi
\in C^2(\mathbb{R})$ is a solution to
%
\begin{eqnarray}
\label{eq:gl104}
\widetilde{\Lambda}(\lambda)&=&\frac{1}{2}\Vert\rho\Vert^2\biggl[\varphi
_{yy}+\frac{1}{1-\lambda}\varphi_y^2\biggr]\nonumber\\[-8pt]\\[-8pt]
&&{}+\inf_{\nu\in\mathbb{R}}\inf_{\eta
\in\Gamma}\{(1-\lambda)l(\eta,\nu,\cdot)+\varphi_yh(\eta,\nu,\cdot)\},\nonumber
\end{eqnarray}
which fulfills the following regularity conditions:
\begin{longlist}[(b)]
\item[(a)] Either the first derivative $\varphi_y$ is bounded or
$\varphi$ is bounded below, and its derivative $\varphi_y$ has at most
linear growth, that is,
\[
\vert\varphi_y(y)\vert\leq C_1(1+\vert y\vert) \qquad\mbox{for some
constant } C_1>0.
\]
\item[(b)]
There exist $C_2,C_3>0$ such that $y\kappa(\eta,y)\leq-C_2y^2+C_3$, where
%
\begin{eqnarray}
\kappa(\eta,y)&:=&g(y)+\frac{\lambda}{1-\lambda}\rho_1\bigl(\theta(y)+\eta
^{11}y+\eta^{21}\bigr)\nonumber\\[-8pt]\\[-8pt]
&&{}+(\rho,\eta^{1\cdot}y+\eta^{2\cdot})+\biggl[\frac
{1}{1-\lambda}\rho_1^2+\rho_2^2\biggr]\varphi_y(y).
\nonumber
\end{eqnarray}
\end{longlist}
\end{assumption}

In full generality, we are unfortunately not able to clarify whether
the EBE~(\ref{eq:gl104}) has such a solution $(\widetilde{\Lambda
}(\lambda),\varphi)$. In Section \ref{sec:existenceergodicBellman}
we are going to state sufficient (but rather restrictive) conditions
under which the existence of a~solution to our EBE (\ref{eq:gl104}) is
already known. Moreover, Section \ref{sec:explicitresults} contains two
case studies with linear drift coefficients, where the solution can be
derived even explicitly. But as illustrated in Section \ref
{sec:geometricOUutility} in case of the geometric OU model, there may
exist multiple such pairs $(\widetilde{\Lambda}(\lambda),\varphi)$,
even beyond the fact that $\varphi$ is determined only except for an
additive constant. However, the verification theorems will require a
certain ``uniform ergodicity condition'' such as
Assumption \ref{th:assump3}(b) for the diffusion $Y$, and this
condition selects the ``good candidate'' for the
optimal growth rate $\Lambda(\lambda)$ (cf. Remark \ref{th:irrelevant}).
\begin{theorem}
\label{th:verification1}
If Assumption \ref{th:assump3} is satisfied, then we get the identity
%
\begin{equation}
\label{eq:gl40}
\widetilde{\Lambda}(\lambda)=\lim_{T\uparrow\infty}\frac{1}{T}\ln\Bigl(\inf
_{\nu\in\mathcal{M}}\inf_{\eta\in\mathcal{C}}V(\eta,\nu
,y_0,T)^{1-\lambda}\Bigr) \qquad\mbox{for any } Y_0=y_0.
\end{equation}
Moreover, the infima at the right-hand side are attained for feedback controls
%
\begin{equation}
\label{eq:gl410}
\eta^*_t:=\eta^*(Y_t), \qquad\nu^*_t:=\nu^*(Y_t), \qquad t\geq0,
\end{equation}
defined in terms of a measurable $\Gamma$-valued function $\eta^*$ and
the function
\[
\nu^*(y):=\nu^*(\eta^*(y),y)=-\eta^{12,*}(y)y-\eta^{22,*}(y)-\rho
_2\varphi_y(y)
\]
such that the infima in (\ref{eq:gl104}) are attained. Thus,
%
\begin{equation}
\label{eq:gl403}
\widetilde{\Lambda}(\lambda)=\lim_{T\uparrow\infty} \frac
{1}{T}\ln(V(\eta^*,\nu^*,y_0,T)^{1-\lambda}).
\end{equation}
In particular, the duality relations for robust utility maximization
yield that
%
\begin{equation}
\label{eq:gl402}\quad
\widetilde{\Lambda}(\lambda)=\lim_{T\uparrow\infty}\frac{1}{T}\ln
U_T(x_0)=\lim_{T\uparrow\infty}\frac{1}{T}\ln U^{Q^{\eta
^*}}_T(x_0) \qquad\mbox{for any } X^{\pi}_0=x_0.
\end{equation}
\end{theorem}
\begin{remark}
\label{th:worstcasemartingale1}
In view of (\ref{eq:gl402}), $Q^{\eta^*}$ can be seen as the \textit
{asymptotic worst-case measure} for robust expected power utility with
parameter $\lambda\in(0,1)$. On the other hand, the probability measure
$P^{\nu^*}$ on $(\Omega,\mathcal{F})$ with Radon--Nikod\'{y}m density
process $(Z_t^{\nu^*})_{t\geq0}$ is a martingale measure which is
equivalent to $Q_0$ on each $\sigma$-algebra $\mathcal{F}_t$, $t>0$. In
view of (\ref{eq:gl403}) and the duality relation (\ref{eq:gl506}) it
can be interpreted as the \textit{asymptotic worst-case martingale measure}.
\end{remark}
\begin{pf*}{Proof of Theorem \ref{th:verification1}}
(1) In order to show that the constant $\widetilde{\Lambda}(\lambda)$
given by the specific solution $(\widetilde{\Lambda}(\lambda),\varphi)$
to the EBE (\ref{eq:gl104}) coincides with the exponential growth rate
of the maximal robust power utility, we first prove that $\widetilde
{\Lambda}(\lambda)$ provides a lower bound for the growth rate. To this
end, we use the duality relation
\[
U_T(x_0)=\frac{1}{\lambda}x_0^\lambda\inf_{\nu\in\mathcal{M}}\inf
_{\eta\in\mathcal{C}} V(\eta,\nu,y,T)^{1-\lambda}
\]
[cf. (\ref{eq:gl47})] with $V$ introduced in (\ref{eq:gl37}), derive
suitable lower bounds for any fixed horizon $T$ and then pass to the limit.

Let $\eta\in\mathcal{C}$, $\nu\in\mathcal{M}$ be fixed controls, and
let $T$ be a given maturity. Then $\tau_n:=\inf\{t\geq0| \vert Y_t\vert
\geq n \mbox{ or } \int
_0^t\nu^2_s\dd s\geq n\}\wedge T$, $n\in\mathbb{N}$, is a localizing
sequence for the local $Q_0$-martingale $(Z_t^{\nu})_{t\leq T}$. This
will allow us to apply the change of measure (\ref{eq:gl138}) locally
and to use the localization Lemma \ref{th:localization1} for $\tau
_n\uparrow T$. In analogy to (\ref{eq:gl507}) we obtain
\[
V(\eta,\nu,y_0,\tau_n)=E_{Q_0}\bigl[\mathcal{E}^{\eta,\nu}_{\tau_n} e^{\int
_0^{\tau_n}l(\eta_t,\nu_t,Y_t)\dd t}\bigr],\qquad n\in\mathbb{N},
\]
where $l$ is the auxiliary function defined in (\ref{eq:gl135}), and where
\[
\mathcal{E}^{\eta,\nu}_{\tau_n}\!=\!\mathcal{E}\biggl(\frac{1}{1\,{-}\,\lambda}\biggl(\int
_0^{\cdot}\!\lambda\theta(Y_u)\,{+}\,\eta^{11}_uY_u\,{+}\,\eta^{21}_u\dd W^1_u\!+\!\int
_0^{\cdot}\!\lambda\nu_u\,{+}\,\eta^{12}_uY_u\,{+}\,\eta^{22}_u\dd W^2_u\biggr)\!\biggr)_{\tau_n}.
\]
To eliminate the It\^{o} exponential $\mathcal{E}^{\eta,\nu}_{\tau_n}$,
we pass to the new probability measure~$R_n^{\eta,\nu}$ on $(\Omega
,\mathcal{F}_T)$ with density process $dR_n^{\eta,\nu}/dQ_0|_{\mathcal
{F}_t}:=\mathcal{E}^{\eta,\nu}_{t\wedge\tau_n}$, $t\in[0,T]$. It
remains to justify this change of measure. For this purpose, note that
the process $\eta\in\mathcal{C}$ takes its values in a compact subset
$\Gamma\subset\mathbb{R}^4$ and that $\theta^2(y)\leq(a_3\vert y\vert
+a_4)^2\leq2(a_3^2y^2+a^2_4)$, due to Assumption \ref
{th:assump1}. Using the definition of~$\tau_n$ we can verify the
Novikov condition (see, e.g., \cite{liptsershiryaev}, Theorem~6.1 and the note after it). This allows us to write
%
\begin{equation}
\label{eq:gl273}
V(\eta,\nu,y_0,\tau_n)=E_{R_n^{\eta,\nu}}\bigl[e^{\int_0^{\tau_n}l(\eta_t,\nu
_t,Y_t)\dd t}\bigr],\qquad n\in\mathbb{N}.
\end{equation}
By Girsanov's theorem, the dynamics of $Y$ follow under $R_n^{\eta
,\nu}$ the SDE
%
\begin{equation}
\label{eq:gl268}
dY_t=h(\eta_t,\nu_t,Y_t)\dd t+\rho\dd W_t^{\eta,\nu} \qquad\mbox{on } \{t\leq
\tau_n\}
\end{equation}
for the drift function $h$ given by (\ref{eq:gl42}) and for a
two-dimensional $R_n^{\eta,\nu}$-Wiener process $W^{\eta,\nu}$. Note
that (\ref{eq:gl273}) can be viewed as a cost functional of an
``expected exponential of integral criterion'' with
dynamics (\ref{eq:gl268}) (cf. page \pageref{eq:gl335}).

Let us next introduce the auxiliary function $\gamma\geq0$ by
%
\begin{eqnarray}
\label{eq:gl261}
\gamma(\eta,\nu,y)&:=&(1-\lambda)l(\eta,\nu,y)+\varphi_y(y)h(\eta,\nu
,y)\nonumber\\[-8pt]\\[-8pt]
&&{}-\inf_{\nu\in\mathbb{R}}\{(1-\lambda)l(\eta,\nu,y)+\varphi
_y(y)h(\eta,\nu,y)\}.
\nonumber
\end{eqnarray}
Inserting the minimizer $\nu^*(\eta,y)$ introduced in (\ref{eq:gl390}),
we then see that $\gamma$ takes the condensed form
%
\begin{equation}
\label{eq:gl262}
\gamma(\eta,\nu,y)=\frac{1}{2}\frac{\lambda}{1-\lambda}\bigl(\nu-\nu
^*(\eta,y)\bigr)^2.
\end{equation}
Later on this representation of $\gamma$ will be crucial to eliminate
the control $\nu$ in the dynamics of $Y$. In terms of $\gamma$ our EBE
(\ref{eq:gl104}) yields the inequality
%
\begin{eqnarray}
\label{eq:gl270}
\widetilde{\Lambda}(\lambda)&\leq&\frac{1}{2}\Vert\rho\Vert^2\biggl[\varphi
_{yy}+\frac{1}{1-\lambda}\varphi_y^2\biggr]\nonumber\\[-8pt]\\[-8pt]
&&{}+(1-\lambda)l(\eta,\nu,\cdot
)+\varphi_yh(\eta,\nu,\cdot)-\gamma(\eta,\nu,\cdot).\nonumber
\end{eqnarray}
By It\^{o}'s formula applied to $\varphi\in C^2(\mathbb{R})$ and to
the dynamics of $Y$ in (\ref{eq:gl268}), this estimate translates on $\{
u\leq\tau_n\}$ into
%
\begin{eqnarray}
&&\varphi(Y_u)-\varphi(y_0)\nonumber\\
&&\qquad=\int_0^u\varphi_y(Y_t)h(\eta_t,\nu_t,Y_t)+\frac{1}{2}\Vert\rho\Vert
^2\varphi_{yy}(Y_t)\dd t+\int_0^u\varphi_y(Y_t)\rho\dd
W_t^{\eta,\nu}\nonumber\\[-8pt]\\[-8pt]
&&\qquad\geq\int_0^u\widetilde{\Lambda}(\lambda)-\frac{1}{2}\frac{1}{1-\lambda
}\Vert\rho\Vert^2\varphi_y^2(Y_t)-(1-\lambda)l(\eta_t,\nu_t,Y_t)\nonumber\\
&&\qquad\quad{}
+\gamma(\eta_t,\nu_t,Y_t)\dd t+\int_0^u\varphi_y(Y_t)\rho\dd W_t^{\eta,\nu}.\nonumber
\end{eqnarray}
Dividing through $1-\lambda$, rearranging the terms and taking the
exponential on both sides, we thus obtain
from (\ref{eq:gl273}) that
%
\begin{eqnarray}
\label{eq:gl224}
&&V(\eta,\nu,y_0,\tau_n)\nonumber\\
&&\qquad\geq E_{R_n^{\eta,\nu}}\biggl[e^{({1}/({1 - \lambda}))(\widetilde{\Lambda}(
\lambda) \tau_n+\varphi( y_0 )-\varphi
( Y_{\tau_n} )+\int_0^{\tau_n} \gamma( \eta_t , \nu_t , Y_t
)\dd t)}\nonumber\\[-8pt]\\[-8pt]
&&\qquad\quad\hspace*{102pt}{}\times
\mathcal{E} \biggl( \int_0^{\cdot} \frac{\varphi_y(Y_t)}{1-\lambda
}\rho\dd W_t^{\eta,\nu}\biggr)_{\tau_n}\biggr]\nonumber\\
&&\qquad= E_{ \overline{R}{}^{\eta,\nu}_n}\bigl[e^{({1}/({1 - \lambda}))(\widetilde
{\Lambda}( \lambda) \tau_n+\varphi( y_0 )-\varphi( Y_{\tau_n} )+\int
_0^{\tau_n} \gamma( \eta
_t , \nu_t , Y_t )\dd t)}\bigr].\nonumber
\end{eqnarray}
Here the last expectation is taken with respect to the probability
measure~$\overline{R}{}^{\eta,\nu}_n$ on $(\Omega,\mathcal{F}_T)$ with
density process
\[
\frac{d\overline{R}{}^{\eta,\nu}_n}{dR_n^{\eta,\nu}}\bigg|_{\mathcal
{F}_t}:=\mathcal{E}\biggl(\int_0^{\cdot}\frac{\rho\varphi_y(Y_u)}{1-\lambda
}\dd W_u^{\eta,\nu}\biggr)_{t\wedge\tau_n}.
\]
Indeed, since $\varphi_y$ grows at most linearly according to
Assumption \ref{th:assump3}(a),
this change of measure can be justified again by Novikov's
condition (cf., e.g.,~\cite{liptsershiryaev}, Theorem 6.1 and
the note after it). By Girsanov's theorem, the factor process $Y$
evolves under $\overline{R}{}^{\eta,\nu}_n$ according to
\[
dY_t=\biggl[h(\eta_t,\nu_t,Y_t)+\frac{1}{1-\lambda}\Vert\rho\Vert^2\varphi
_y(Y_t)\biggr]\dd t+\rho\dd\overline{W}{}^{\eta,\nu}_t \qquad\mbox{on } \{t\leq\tau
_n\},
\]
where $\overline{W}{}^{\eta,\nu}$ denotes a two-dimensional $\overline
{R}{}^{\eta,\nu}_n$-Wiener process. But these dynamics still depend on
the irrepressible control $\nu$. To eliminate this dependence, we apply
once more a Girsanov transformation. Consider the probability measure
$R_n^{\eta}$ on $(\Omega,\mathcal{F}_T)$ with density process
\[
\frac{d\widehat{R}_n^{\eta}}{d\overline{R}{}^{\eta,\nu}_n}\bigg|_{\mathcal
{F}_t}:=\mathcal{E}\biggl(\int_0^{\cdot}\frac{\lambda}{1-\lambda}\bigl(\nu^*(\eta
_s,Y_s)-\nu_s\bigr)\dd\overline{W}{}^{2,\eta,\nu}_t\biggr)_{t\wedge\tau_n}.
\]
Verifying once more Novikov's condition, we see that $R_n^{\eta}$ is
well defined, and so the inequality (\ref{eq:gl224}) translates into
%
\begin{eqnarray}
\label{eq:gl259}
\quad&&V( \eta, \nu, y_0 , \tau
_n ) \nonumber\\[-8pt]\\[-8pt]
\quad&&\qquad\geq E_{\widehat
{R}_n^{\eta}}\biggl[e^{({1}/({1-\lambda}))(\widetilde{\Lambda}(\lambda) \tau
_n+\varphi(y_0)-\varphi(Y_{\tau_n})+\int_0^{\tau_n} \gamma(\eta_t,\nu
_t,Y_t)\dd t)}\frac{d\overline{R}{}^{\eta,\nu
}_n}{d\widehat{R}_n^{\eta}}\bigg|_{\mathcal{F}_{\tau_n}}
\biggr].\nonumber
\end{eqnarray}
Moreover, Girsanov's theorem yields that the dynamics of $Y$ under
$\widehat{R}_n^{\eta}$ on $\{t\leq\tau_n\}$ takes the form
\[
dY_t=\biggl[h(\eta_t,\nu_t,Y_t)+\frac{1}{1-\lambda}\Vert\rho\Vert^2\varphi
_y(Y_t)+\frac{\lambda}{1-\lambda}\rho_2\bigl(\nu^*(\eta_t,Y_t)-\nu
_t\bigr)\biggr]\dd t+\rho\dd\widehat{W}^{\eta}_t
\]
in terms of the two-dimensional $\widehat{R}_n^\eta$-Wiener process
$\widehat{W}^{\eta}$. Recalling\break from (\ref{eq:gl42}) and (\ref
{eq:gl390}) the definitions of the drift function $h$ and of the
minimi\-zer~$\nu^*(\eta,y)$, a straightforward computation shows that
this SDE is equivalent to\looseness=1
%
\begin{equation}
\label{eq:gl157}
dY_t=\kappa(\eta_t,Y_t)\dd t+\rho\dd\widehat{W}_t^{\eta},
\end{equation}\looseness=0
where\vspace*{1pt} $\kappa$ denotes the auxiliary function introduced in
Assumption \ref{th:assump3}(b). To eliminate the density $d\overline
{R}{}^{\eta,\nu}_n/d\widehat{R}_n^{\eta}|_{\mathcal{F}_{\tau_n}}$, we
define $p:=\frac{\lambda-1}{\lambda}<0$ and apply H\"{o}lder's
inequality with $1/p+1/q=1$ to (\ref{eq:gl259}) (see, e.g.,
\cite{hewittstromberg},\vadjust{\goodbreak} page 191, for an extension of the classical result
to $p<0$, $q\in(0,1)$). This
leads to
%
\begin{eqnarray}
\label{eq:gl260}\quad
V( \eta, \nu, y_0 , \tau
_n )& \geq & E_{\widehat
{R}_n^{\eta}}\bigl[e^{({q}/({1-\lambda}))(\widetilde{\Lambda}(\lambda)\tau
_n+\varphi(y_0)-\varphi(Y_{\tau_n} ))}\bigr]^{1/q}\nonumber\\[-8pt]\\[-8pt]
&&{} \times E_{\widehat{R}_n^{\eta}}\biggl[
\biggl(\frac{d\overline{R}{}^{(\eta,\nu)}_n}{d\widehat{R}_n^{(\eta)}}
\bigg|_{\mathcal{F}_{\tau_n}}e^{({1}/({1-\lambda}))\int_0^{\tau_n}\gamma
(\eta_t,\nu_t,Y_t)\dd t}\biggr)^p\biggr]^{1/p}.\nonumber
\end{eqnarray}
But in view of (\ref{eq:gl262}) and our choice of $p$ we see that
\begin{eqnarray*}
&&\biggl(\frac{d\overline{R}{}^{\eta,\nu}_n}{d\widehat{R}_n^{\eta}}
\bigg|_{\mathcal{F}_{\tau_n}}e^{({1}/({1-\lambda}))\int_0^{\tau_n}\gamma
(\eta_t,\nu_t,Y_t)\dd t}\biggr)^p \\
&&\qquad=\mathcal{E}\biggl( \int_0^{\cdot} \frac{p\lambda
}{1-\lambda}\bigl(\nu
^*( \eta_t , Y_t )-\nu_t\bigr)\dd\widehat{W}_t^{2,\eta}\biggr)_{\tau_n}.
\end{eqnarray*}
Since the It\^{o} exponential of a local martingale is always a
supermartingale, it follows
that the expectation in (\ref{eq:gl260}) is less than $1$. Raised to
the power of $1/p<0$, this estimate is reversed, and we obtain
%
\begin{equation}
\label{eq:gl257}
\quad V(\eta,\nu,y_0,\tau_n)\geq E_{\widehat{R}_n^{\eta}}\bigl[e^{
({q}/({1-\lambda}))(\widetilde{\Lambda}(\lambda)\tau_n+\varphi(y_0)-\varphi
(Y_{\tau_n}))}\bigr]^{1/q},\qquad n\in\mathbb{N}.
\end{equation}
In our next step, we shall extend the measures $\widehat{R}_n^{\eta
}|_{\mathcal{F}_{\tau_n}}$, $n\in\mathbb{N}$, to a probability measure
$\widehat{R}^{\eta}$ on the $\sigma$-field $\mathcal{F}_T$ whose
restrictions to $\mathcal{F}_{\tau_n}$ are equal to $\widehat{R}_n^{\eta
}|_{\mathcal{F}_{\tau_n}}$ for all $n\in\mathbb{N}$. To this end, note
that the sequence $\tau_n$ increases to $T$ and that the family
$(\widehat{R}_n^{\eta}|_{\mathcal{F}_{\tau_n}})_{n\in\mathbb{N}}$ is
consistent in the sense that $\widehat{R}_{n+1}^{\eta}(A)=\widehat
{R}_n^{\eta}(A)$ for all $A\in\mathcal{F}_{\tau_n}$ since
%
\begin{eqnarray}
\label{eq:gl531}
\frac{d\widehat{R}_n^{\eta}}{dQ_0}\bigg|_{\mathcal{F}_{\tau_n}}&=&\frac
{d\widehat{R}_n^{\eta}}{d\overline{R}{}^{\eta,\nu}_n}\bigg|_{\mathcal
{F}_{\tau_n}}\frac{d\overline{R}{}^{\eta,\nu}_n}{dR_n^{\eta,\nu}}
\bigg|_{\mathcal{F}_{\tau_n}}\frac{dR_n^{\eta,\nu}}{dQ_0}\bigg|_{\mathcal
{F}_{\tau_n}}\nonumber\\
&=&\mathcal{E}\biggl(\int_0^{\cdot}\frac{1}{1-\lambda}\bigl(\lambda\theta(Y_u)+\eta
^{11}_uY_u+\eta^{21}_u+\varphi_y(Y_u)\rho_1\bigr)\dd
W^1_u\\
&&\hspace*{80.4pt}{}
+\int_0^{\cdot}\eta^{12}_uY_u+\eta^{22}_u+\varphi_y(Y_u)\rho_2\dd
W^2_u\biggr)_{\tau_n},\nonumber
\end{eqnarray}
$n\in\mathbb{N}$,\vspace*{1pt} is a discrete-time $Q_0$-martingale. Thus the existence of a unique
extension~$\widehat{R}^{\eta}$ to $\sigma(\bigcup_{n\in\mathbb{N}}\mathcal
{F}_{\tau_n})=\mathcal{F}_T$ follows from \cite{parthasarathy},
Theorem V.4.2. More directly, (\ref{eq:gl531}) suggests that we should
define the probability measure $\widehat{R}^{\eta}$ on~$(\Omega,\mathcal
{F}_T)$ by
%
\begin{eqnarray}
\label{eq:gl424}\quad
\frac{d\widehat{R}^{\eta}}{dQ_0}\bigg|_{\mathcal{F}_T}&:=&
\mathcal{E}\biggl(\int_0^{\cdot}\frac{1}{1-\lambda}\bigl(\lambda\theta(Y_u)+\eta
^{11}_uY_u+\eta^{21}_u+\varphi_y(Y_u)\rho_1\bigr)\dd
W^1_u\nonumber\\[-8pt]\\[-8pt]
&&\hspace*{80.4pt}{} +\int_0^{\cdot}\eta^{12}_uY_u+\eta^{22}_u+\varphi_y(Y_u)\rho_2\dd
W^2_u\biggr)_{T}.\nonumber
\end{eqnarray}
Since the functions $\theta$, $\varphi$ grow, at most, linearly, it
follows similarly to page~\pageref{eq:gl381} that $\widehat{R}^{\eta}$
is well defined, that is, $E_{Q_0}[d\widehat{R}^{\eta} / dQ_0|_{\mathcal
{F}_T}]=1$. In particular, the
corresponding It\^{o} exponential is a $Q_0$-martingale up to time $T$,
and in view of~(\ref{eq:gl531}) this yields $\widehat{R}^{\eta
}|_{\mathcal{F}_{\tau_n}}=\widehat{R}_n^{\eta}|_{\mathcal{F}_{\tau_n}}$
for all $n\in\mathbb{N}$. We thus see that estimate~(\ref{eq:gl257}) is
equivalent to
\[
V(\eta,\nu,y_0,\tau_n)\geq E_{\widehat{R}^{\eta}}\bigl[e^{
({q}/({1-\lambda}))(\widetilde{\Lambda}(\lambda)\tau_n+\varphi(y_0)-\varphi
(Y_{\tau_n}))}\bigr]^{1/q} \qquad\mbox{for any } n\in\mathbb{N}.
\]
Now we are ready to replace the stopping times $\tau_n$ by the
deterministic time $T$ by passing to the limit $n\uparrow\infty$.
Indeed, as shown in Lemma \ref{th:localization1}, the left-hand
side increases to $V(\eta,\nu,y_0,T)$ as $n\uparrow\infty$ (cf.
Lemma \ref{th:localization1}). Applying Fatou's lemma and then
Jensen's inequality to the rightmost expectation, we now obtain the
lower bound
\begin{eqnarray*}
V(\eta,\nu, y_0 , T ) &\geq&
E_{\widehat{R}^{\eta}}\bigl[e^{({1}/({1-\lambda
}))(\widetilde{\Lambda}(\lambda)T+\varphi(y_0)-\varphi(Y_T))}\bigr]\\
&\geq&
e^{({1}/({1-\lambda}))(\widetilde{\Lambda
}(\lambda)T+\varphi(y_0)+E_{\widehat{R}^{\eta}}[-\varphi(Y_T)])}
\end{eqnarray*}
for any finite horizon $T$ and for all controls $\eta\in\mathcal{C}$,
$\nu\in\mathcal{M}$. Taking the scaling $\frac{1}{T}\ln(\cdot
)^{1-\lambda}$ on both sides and passing to the limit $T\uparrow\infty
$, this yields
\[
\mathop{\underline{\lim}}_{T\uparrow\infty} \frac{1}{T}\ln
\Bigl(\inf_{\nu\in\mathcal{M}}\inf_{\eta\in\mathcal{C}}V(\eta,\nu
,y_0,T)^{1-\lambda}\Bigr)\geq\widetilde{\Lambda}(\lambda)+\mathop
{\underline{\lim}}_{T\uparrow\infty} \frac{1}{T}\inf_{\eta\in
\mathcal{C}}E_{\widehat{R}^{\eta}}[-\varphi(Y_T)].
\]
Thus the constant $\widetilde{\Lambda}(\lambda)$ provides a lower bound if
%
\begin{equation}
\label{eq:gl263}
\lim_{T\uparrow\infty} \frac{1}{T}\inf_{\eta\in\mathcal
{C}}E_{\widehat{R}^{\eta}}[-\varphi(Y_T)]=0.
\end{equation}
Indeed, Assumption \ref{th:assump3}(a) ensures that $\varphi$
grows at most quadratically, that is, there exists some constant
$K_1>0$ with $\vert\varphi(y)\vert\leq K_1(1+y^2)$. Therefore, we have
the bounds
%
\begin{equation}
\label{eq:gl264}\qquad
-K_1\Bigl(1+\sup_{\eta\in\mathcal{C}}E_{\widehat{R}^{\eta}}[Y_T^2]\Bigr)\leq\inf
_{\eta\in\mathcal{C}}E_{\widehat{R}^{\eta}}[-\varphi(Y_T)]\leq
K_1\Bigl(1+\sup_{\eta\in\mathcal{C}}E_{\widehat{R}^{\eta}}[Y_T^2]\Bigr).
\end{equation}
Recall now from (\ref{eq:gl157}) that $Y$ evolves under $\widehat
{R}^{\eta}$, $\eta\in\mathcal{C}$, according to the SDE
\[
dY_t=\kappa(\eta_t,Y_t)\dd t+\rho(Y_t)\dd\widehat{W}_t^{\eta}.
\]
Due to Assumption \ref{th:assump3}(b) there exist constants
$C_2,C_3>0$ such that the drift function $\kappa$ satisfies $y\kappa
(\eta,y)\leq-C_2y^2+C_3$ for all $\eta\in\Gamma$. Therefore,
Lemma~\ref{th:boundedmoments} ensures that
\[
\sup_{T\geq0}\sup_{\eta\in\mathcal{C}}E_{\widehat{R}^{\eta
}}[Y_T^2]\leq y_0^2+\mbox{const.}<\infty.
\]
But in view of (\ref{eq:gl264}) this implies (\ref{eq:gl263}), and hence
%
\begin{equation}
\label{eq:gl265}
\mathop{\underline{\lim}}_{T\uparrow\infty} \frac{1}{T}\ln
\Bigl(\inf_{\nu\in\mathcal{M}}\inf_{\eta\in\mathcal{C}}V(\eta,\nu
,y_0,T)^{1-\lambda}\Bigr)\geq\widetilde{\Lambda}(\lambda).
\end{equation}

(2) In the second part we identify controls $\eta^*\in\mathcal{C}$ and
$\nu^*\in\mathcal{M}$ such that
%
\begin{equation}
\label{eq:gl425}
\widetilde{\Lambda}(\lambda)=\lim_{T\uparrow\infty}\frac{1}{T}\ln
(V(\eta^*,\nu^*,y_0,T)^{1-\lambda}).
\end{equation}
Together with (\ref{eq:gl265}) this implies (\ref{eq:gl40}). Indeed, by
compactness of $\Gamma$ and continuity of the functions $l$, $h$ and
$\nu^*(\cdot,y)$
with respect to $\eta$, there exists
%
\begin{equation}
\eta^*(y)\in\mathop{\arg\min}_{\eta\in\Gamma} \{(1-\lambda
)l(\eta,\nu^*(y,\eta),y)+\varphi_y(y)h(\eta,\nu^*(y,\eta),y)\}.
\end{equation}
By a measurable selection argument $\eta^*(\cdot)$ can be chosen as a
measurable function. Set $\nu^*(y):=\nu^*(\eta^*(y),y)$ [cf. (\ref
{eq:gl390})], and let $\eta^*$, $\nu^*$ be the feedback controls defined
by $\eta^*_t:=\eta^*(Y_t)$, $\nu_t^*:=\nu^*(Y_t)$, $t\geq0$. In that
case, we have $\eta^*\in\mathcal{C}$, and one easily proves that the
process $\nu^*$ belongs to the class $\mathcal{M}$.

In order to verify (\ref{eq:gl425}), we now proceed as in part (1). As
in (\ref{eq:gl273}) we obtain
\[
V(\eta^*,\nu^*,y_0,T)=E_{R^{\eta^*,\nu^*}}\bigl[e^{\int_0^T l(\eta^*_t,\nu
^*_t,Y_t)\dd t}\bigr].
\]
The measure $R^{\eta^*,\nu^*}$ is defined on $(\Omega,\mathcal{F}_T)$
in terms of the density $\mathcal{E}_T^{\eta^*,\nu^*}$. Since $\nu
^*(\eta,\cdot)$ grows at most linearly, this change of measure can be
justified in analogy to page \pageref{eq:gl381}. By Girsanov's
theorem, the dynamics of $Y$ under $R^{\eta^*,\nu^*}$ follow the SDE
%
\begin{equation}
\label{eq:gl279}
dY_t=h(\eta^*_t,\nu^*_t,Y_t)\dd t+\rho\dd W_t^{\eta^*,\nu^*},
\end{equation}
where the drift\vspace*{1pt} function $h$ is given by (\ref{eq:gl42}), and where
$(W_t^{\eta^*,\nu^*})_{t\leq T}$ is a two-dimensional Wiener process
under $R^{\eta^*,\nu^*}$ (cf. page \pageref{eq:gl268}). Using the
specific controls $\eta^*$, $\nu^*$, the auxiliary function $\gamma$ in
(\ref{eq:gl261}) satisfies $\gamma(\eta^*_t,\nu^*_t,Y_t)=0$, and we
also obtain equality in (\ref{eq:gl270}). Along the lines of part (1)
this implies
\begin{eqnarray*}
&&V( \eta^* , \nu^* , y_0 , T ) \\
&&\qquad= E_{R^{\eta
^*,\nu^*}}\biggl[e^{({1}/({1-\lambda}))(\widetilde{\Lambda}(\lambda)
T+\varphi(y_0)-\varphi(Y_T))}\mathcal{E}\biggl(\int_0^{\cdot}\frac{\varphi
_y(Y_t)}{1-\lambda}\rho\dd W_t^{\eta^*,\nu^*}\biggr)_T\biggr]
\end{eqnarray*}
in analogy to (\ref{eq:gl224}). Once more the It\^{o} exponential is
interpreted as the density of a new probability measure $\widehat
{R}^{\eta^*}$ on $(\Omega,\mathcal{F}_T)$. Since the drift function
$h(\eta^*(\cdot),\nu^*(\cdot),\cdot)$ of $Y$ under $R^{\eta^*,\nu^*}$
only depends on the control $\eta^*$ and satisfies the linear growth
condition $\vert h(\eta^*(y),\nu^*(y),y)\vert\leq K_2(1+\vert y\vert)$,
we may proceed in analogy to page~\pageref{eq:gl381} to justify this
change of measure. Then we get
%
\begin{equation}
\label{eq:gl280}
V(\eta^*,\nu^*,y_0,T)=e^{({1}/({1-\lambda}))(\widetilde{\Lambda
}(\lambda) T+\varphi(y_0)}E_{\widehat{R}^{\eta^*}}\bigl[e^{-
({1}/({1-\lambda}))\varphi(Y_T)}\bigr].
\end{equation}
Moreover, by Girsanov's theorem, the dynamics of $Y$ with respect to
$\widehat{R}^{\eta^*}$ are given by
%
\begin{equation}
\label{eq:gl283}
dY_t=\kappa(\eta^*_t,Y_t)\dd t+\rho\dd\widehat{W}_t^{\eta^*},
\end{equation}
where $(\widehat{W}_t^{\eta^*})_{t\leq T}$ is a two-dimensional Wiener
process, and where the drift function $\kappa$ satisfies
Assumption \ref{th:assump3}(b). In analogy to part (1), we now take
the scaling $\frac{1}{T}\ln(\cdot)^{1-\lambda}$ on both sides of (\ref
{eq:gl280}) and then pass to the limit $T\uparrow\infty$. For this
purpose, note that
\[
\sup_{T\geq0}E_{\widehat{R}^{\eta^*}}[Y^2_T]<\infty\quad\mbox{and
that}\quad
\sup_{T\geq0}E_{\widehat{R}^{\eta^*}}[\exp(k\vert Y_T\vert)]<\infty
\qquad\mbox{for any } k\in\mathbb{R},\vspace*{-2pt}
\]
due to Assumption \ref{th:assump3}(b) and Lemma \ref
{th:boundedmoments} applied to the SDE (\ref{eq:gl283}).
If $\varphi_y$ is bounded and consequently $\vert\varphi(y)\vert\leq
K_3(1+\vert y\vert)$, then this implies the uniform upper bound
\[
\sup_{T\geq0}E_{\widehat{R}^{\eta^*}}\biggl[\exp\biggl(-\frac{1}{1-\lambda
}\varphi(Y_T)\biggr)\biggr]\leq\sup_{T\geq0}E_{\widehat{R}^{\eta^*}}\biggl[\exp\biggl(\frac
{1}{1-\lambda}K_3(1+\vert Y_T\vert)\biggr)\biggr]<\infty.\vspace*{-2pt}
\]
This uniform boundedness among all $T$ clearly also holds, if $\varphi$
is bounded below. In particular, the identity (\ref{eq:gl280})
translates into
\[
\lim_{T\uparrow\infty}\frac{1}{T}\ln( V(\eta^*,\nu^*,y_0,T)^{1-\lambda
})=\widetilde{\Lambda}(\lambda).\vspace*{-2pt}
\]
Thus we have shown (\ref{eq:gl425}). This ends the proof of (\ref{eq:gl40}).

(3) In our last step we return to the initial problem of robust utility
maximization. The finite horizon duality relation (\ref{eq:gl506})
holds for any (regular) convex class of measures, and in particular for
the one-point set $\{Q^{\eta^*}\}$. In analogy to (\ref{eq:gl47}) it
thus follows that the maximal value for expected power utility in the
specific model $Q^{\eta^*}$ satisfies the duality formula
\[
U^{Q^{\eta^*}}_T(x_0)=\frac{1}{\lambda}x_0^\lambda
\Bigl(\inf_{\nu\in\mathcal{M}}V(\eta^*,\nu,y_0,T)\Bigr)^{1-\lambda}.\vspace*{-2pt}
\]
Using this representation and the duality relation (\ref{eq:gl47}) for
the whole set $\mathcal{Q}$, we obtain (\ref{eq:gl402}) immediately
from (\ref{eq:gl40}) and (\ref{eq:gl403}).\vspace*{-2pt}
\end{pf*}

Theorem \ref{th:verification1} shows that the solution
$(\widetilde{\Lambda}(\lambda),\varphi)$ to the EBE (\ref{eq:gl104})
specified in Assumption \ref{th:assump3} describes the
exponential growth of the maximal robust power utility $U_T(x_0)$ as
$T\uparrow\infty$. We have also seen that the maximal utility in the
specific model $Q^{\eta^*}$ grows at the same rate as $U_T(x_0)$. In
the next step we shall use these facts in order to identify an optimal
long-term investment strategy $\pi^*\in\mathcal{A}$. For this purpose,
we introduce the additional regularity.\vspace*{-2pt}

\begin{assumption}
\label{th:assump4}
Let $(\widetilde{\Lambda}(\lambda),\varphi)$ be the solution to the
EBE (\ref{eq:gl104}) introduced in Assumption \ref{th:assump3},
and let $\eta^*$ denote the corresponding minimizing function. Then the
function $\widetilde{\kappa}$ defined by
%
\begin{eqnarray}
\widetilde{\kappa}(\eta,y)&:=&g(y)+\frac{\lambda}{1-\lambda}\rho_1\bigl(\theta
(y)+\eta^{11,*}(y)y+\eta^{21,*}(y)\bigr)\nonumber\\[-10pt]\\[-10pt]
&&{}+(\rho,\eta^{1\cdot}y+\eta^{2\cdot})+\biggl[\frac{1}{1-\lambda}\rho_1^2+\rho
_2^2\biggr]\varphi_y(y)
\nonumber\vspace*{-2pt}
\end{eqnarray}
satisfies $y\widetilde{\kappa}(\eta,y)\leq-C_4y^2+C_5$ for all $\eta\in
\Gamma$ with constants $C_4,C_5>0$.\vadjust{\goodbreak}
\end{assumption}
\begin{theorem}
\label{th:identification}
Under the regularity Assumptions \ref{th:assump3} and \ref
{th:assump4} we have:
\begin{longlist}
\item The value $\widetilde{\Lambda}(\lambda)$ given by the
solution to the EBE (\ref{eq:gl104}) can be identified as the optimal
exponential growth rate
\[
\Lambda(\lambda)=\sup_{\pi\in\mathcal{A}}\mathop{\overline
{\lim}}_{T\uparrow\infty} \frac{1}{T}\ln\inf_{Q^{\eta}\in\mathcal
{Q}} E_{Q^{\eta}}[(X_T^{\pi})^{\lambda}]
\]
for robust expected power utility. In particular, (\ref{eq:gl402}) implies
\[
\Lambda(\lambda)=\lim_{T\uparrow\infty}\frac{1}{T}\ln U_T(x_0)=\lim
_{T\uparrow\infty}\frac{1}{T}\ln U_T^{Q^{\eta^*}}(x_0),
\]
where $Q^{\eta^*}\in\mathcal{Q}$ is defined in terms of the control
$\eta^*$ in (\ref{eq:gl410}).
\item In the specific model $Q^{\eta^*}$, the maximal
growth rate of power utility
\[
\Lambda_{Q^{\eta^*}}(\lambda):= \sup_{\pi\in\mathcal{A}}
\mathop{\overline{\lim}}_{T\uparrow\infty} \frac{1}{T}\ln E_{Q^{\eta
^*}}[(X_T^{\pi})^{\lambda}]
\]
coincides with $\Lambda(\lambda)$.
\item Let $\pi_t^*=\pi^*(Y_t)$, $t\geq0$, be the trading
strategy defined in terms of the function (\ref{eq:gl24}). Then $\pi^*$
belongs to class $\mathcal{A}$, and it satisfies the optimality condition
%
\begin{equation}
\label{eq:gl71}
\Lambda(\lambda)=\lim_{T\uparrow\infty} \frac{1}{T}\ln\inf
_{Q^{\eta}\in\mathcal{Q}} E_{Q^{\eta}}[(X_T^{\pi^*})^\lambda]=\lim
_{T\uparrow\infty}\frac{1}{T}\ln E_{Q^{\eta^*}}[(X_T^{\pi^*})^\lambda].
\end{equation}
\end{longlist}
In other words, the strategy $\pi^*$ and the measure $Q^{\eta^*}\in
\mathcal{Q}$ form a saddle point for the robust optimization problem
(\ref{eq:gl27}).
\end{theorem}
\begin{pf} (1) Theorem \ref{th:verification1} shows that the
maximal power utility $U_T^{Q^{\eta^*}}(x_0)$ in the
specific model $Q^{\eta^*}$ grows exponentially with rate $\widetilde
{\Lambda}(\lambda)$, that is,
\[
\widetilde{\Lambda}(\lambda)=\lim_{T\uparrow\infty}\frac{1}{T}\ln
U_T^{Q^{\eta^*}}(x_0)=\lim_{T\uparrow\infty}\frac
{1}{T}\ln\sup_{\pi\in\mathcal{A}_T}E_{Q^{\eta^*}}[(X_T^{\pi})^\lambda].
\]
Since $\mathcal{A}\subseteq\mathcal{A}_T$, this implies
\[
\widetilde{\Lambda}(\lambda)\geq\sup_{\pi\in\mathcal{A}}\mathop
{\overline{\lim}}_{T\uparrow\infty} \frac{1}{T}\ln E_{Q^{\eta
^*}}[(X_T^{\pi})^{\lambda}]\geq\sup_{\pi\in\mathcal{A}}\mathop
{\overline{\lim}}_{T\uparrow\infty} \frac{1}{T}\ln\inf_{Q^{\eta
}\in\mathcal{Q}}E_{Q^{\eta}}[(X_T^{\pi})^{\lambda}]=\Lambda(\lambda).
\]
In order to verify that this chain of inequalities is indeed a series
of equalities, it suffices to show that $\pi^*$ belongs to $\mathcal
{A}$, and that
%
\begin{equation}
\label{eq:gl111}
\widetilde{\Lambda}(\lambda)\leq\mathop{\underline{\lim}}_{T\uparrow\infty} \frac{1}{T}\ln\inf_{Q^{\eta}\in\mathcal
{Q}} E_{Q^{\eta}}[(X_T^{\pi^*})^\lambda].
\end{equation}
This yields the converse inequality $\widetilde{\Lambda}(\lambda)\leq
\Lambda(\lambda)$, and hence the identity $\widetilde{\Lambda}(\lambda
)=\Lambda(\lambda)=\Lambda_{Q^{\eta^*}}(\lambda)$.
In particular, the strategy $\pi^*$ satisfies
(\ref{eq:gl71}).\vspace*{2pt}

Let us first show that $\pi^*$ is admissible in the sense of
Definition \ref{th:definitionstrategyadmissible}. For this purpose,
note that the adapted process $\pi^*_t=\pi^*(Y_t)$, $t\geq0$, admits
continuous paths and that the unique strong solution to (\ref{eq:gl89})
takes the form
%
\begin{equation}
\label{eq:finalformula1}\qquad
X_t^{\pi^*}=x_0e^{\int_0^t\pi^*_u\sigma\dd W_u^{1,\eta}+\int
_0^Tr(Y_u)+\sigma\pi^*_u[\theta(Y_u)+\eta_u^{11}Y_u+\eta^{21}_u]-
({1/2})\sigma^2(\pi^*_u)^2\dd u)}>0
\end{equation}
for any $t\geq0$. Thus the processes defined by the number of shares,
\[
\xi^{*,0}_t=\frac{X_t^{\pi^*}(1-\pi^*_t)}{S^0_t} \quad\mbox{and}\quad
\xi^{*,1}_t=\frac{X_t^{\pi^*}\pi^*_t}{S^1_t},\qquad t\geq0,
\]
are continuous and\vspace*{1pt} adapted to the Brownian filtration, hence
predictable. Moreover, the integrals in (\ref{eq:gl400}) are
well defined for $\xi^*=(\xi^{*,0},\xi^{*,1})$. In other words,
$\pi^*$~associated with $\xi^*$ is an admissible long-term investment process.

To verify (\ref{eq:gl111}), we derive suitable lower bounds for $\inf
_{Q^{\eta}\in\mathcal{Q}} E_{Q^{\eta}}[(X_T^{\pi^*})^\lambda]$ for any
finite horizon $T$ and then pass to the limit. We first argue for a~fixed control $\eta\in\mathcal{C}$ and the corresponding model $Q^{\eta
}\in\mathcal{Q}$.
Representation~(\ref{eq:finalformula1}) yields the decomposition
%
\begin{equation}
\label{eq:gl413}
E_{Q^{\eta}}[(X_T^{\pi^*})^\lambda]=x_0^\lambda E_{Q^{\eta}}\biggl[\mathcal
{E}\biggl(\int_0^{\cdot}\lambda\sigma\pi^*_t\dd W_t^{1,\eta}\biggr)_Te^{\int_0^T
\widetilde{l}(\pi^*_t,\eta_t,Y_t)\dd t}\biggr],
\end{equation}
where we use, as in (\ref{eq:gl108}), the function $\widetilde{l}$. In
order to eliminate the It\^{o} exponential, we introduce a new
probability measure $\overline{Q}{}^{\eta}$ on $(\Omega,\mathcal{F}_T)$
with density
%
\begin{equation}
\label{eq:gl266}
\frac{d\overline{Q}{}^{\eta}}{dQ^{\eta}}\bigg|_{\mathcal{F}_T}:=\mathcal
{E}\biggl(\int_0^{\cdot}\lambda\sigma\pi^*_t\dd W_t^{1,\eta}\biggr)_T=\mathcal
{E}\biggl(\int_0^{\cdot}\lambda\sigma\pi^*(Y_t)\dd W_t^{1,\eta}\biggr)_T.
\end{equation}
This requires us to verify $E_{Q^\eta}[d\overline{Q}{}^{\eta}/dQ^{\eta
}|_{\mathcal{F}_T}]=1$. Indeed, the factor process~$Y$ evolves under
$Q^{\eta}$ according to the SDE (\ref{eq:gl267}), and the drift
function satisfies
\[
\vert g(y)+(\rho,\eta^{1\cdot}y+\eta^{2\cdot})\vert^2\leq K_1(1+y^2),
\]
due to Assumption \ref{th:assump1} and compactness of $\Gamma
\subset\mathbb{R}^4$. Thus, by Lemma \ref{th:liptserersatz},
there exists some constant $K_2>0$ such that $\sup_{0\leq t\leq
T}E_{Q^{\eta}}[\exp(K_2 Y_t^2)]<\infty$.
Since $ \vert\pi^*(y)\vert\leq K_3(1+\vert y\vert)$, this implies $\sup
_{ t\leq T}E_{Q^{\eta}}[\exp(\delta(\lambda\sigma\pi^*(Y_t))^2]<\infty
$ as soon as $\delta>0$ is chosen sufficiently small. Therefore,
\cite{liptsershiryaev}, Example 3 of Section~6.2, guarantees that (\ref
{eq:gl266}) defines a probability measure on $(\Omega,\mathcal{F}_T)$.
In particular, equation (\ref{eq:gl413}) becomes equivalent to
%
\begin{equation}
\label{eq:gl72}
E_{Q^{\eta}}[(X_T^{\pi^*})^\lambda]=x_0^\lambda E_{\overline{Q}{}^{\eta
}}\bigl[e^{\int_0^T \widetilde{l}(\pi^*_t,\eta_t,Y_t)\dd t}\bigr].
\end{equation}
By Girsanov's theorem, the factor process $Y$ follows under
$\overline{Q}{}^{\eta}$ the SDE
%
\begin{equation}
\label{eq:gl109}
dY_t=\widetilde{h}(\pi^*_t,\eta_t,Y_t)\dd t+\rho\dd\overline{W}{}^{\eta
}_t, \qquad t\leq T,\qquad Y_0=y_0.
\end{equation}
Here $(\overline{W}{}^{\eta}_t)_{t\leq T}$ is a two-dimensional $\overline
{Q}{}^{\eta}$-Wiener process and the drift function~$\widetilde{h}$ is
defined by (\ref{eq:gl98}). Note that the right-hand side of (\ref
{eq:gl72}) can be viewed as a cost functional of an ``exponential
of integral criterion'' with
dynamics~(\ref{eq:gl109}).\vadjust{\goodbreak}

In terms of the functions $\widetilde{l}$ and $\widetilde{h}$ the EBE
(\ref{eq:gl104}) for the pair $(\widetilde{\Lambda}(\lambda),\varphi)$
can be rewritten as
%
\begin{equation}
\label{eq:gl106}
\widetilde{\Lambda}(\lambda)=\frac{1}{2}\Vert\rho\Vert^2[\varphi
_{yy}+\varphi^2_y]+\inf_{\eta\in\Gamma}\{\widetilde{l}(\pi^*,\eta,\cdot
)+\varphi_y\widetilde{h}(\pi^*,\eta,\cdot)\}.
\end{equation}
For clarity of exposition the precise arguments are postponed to part
(2) of this proof. We now proceed in analogy\vspace*{2pt} to the proof of
Theorem \ref{th:verification1}. Note that the roles played by $l$, $h$
are taken over by $\widetilde{l}$, $\widetilde{h}$.

Applying It\^{o}'s formula to $\varphi\in C^2(\mathbb{R})$ and to
the dynamics (\ref{eq:gl109}) we obtain
\[
\varphi( Y_T ) = \varphi( y_0 ) + \int_0^T \varphi_y( Y_t )\widetilde
{h}( \pi
^*_t , \eta_t , Y_t ) + \frac
{1}{2}\Vert\rho\Vert^2\varphi_{yy}( Y_t )\dd t + \int_0^T
\varphi_y( Y_t )\rho\dd\overline
{W}{}^{\eta}_t.
\]
The alternative version (\ref{eq:gl106}) of our EBE thus yields the inequality
\[
\varphi(Y_T)\geq\varphi(y_0)+\int_0^T\widetilde{\Lambda}(\lambda
)-\widetilde{l}(\pi^*_t,\eta_t,Y_t)\dd t+\ln\mathcal{E}\biggl(\int_0^{\cdot
}\varphi_y(Y_t)\rho\dd\overline{W}{}^{\eta}_t\biggr)_T.
\]
Rearranging the terms and taking the exponential on both sides, (\ref
{eq:gl72}) allows us to deduce that
\begin{eqnarray*}
E_{Q^{\eta}}[(X_T^{\pi^*})^\lambda]&=&x_0^\lambda E_{\overline{Q}{}^{\eta
}}\bigl[e^{\int_0^T \widetilde{l}(\pi^*_t,\eta_t,Y_t)\dd t}\bigr]\\
&\geq&
x_0^\lambda e^{\widetilde{\Lambda}(\lambda) T+\varphi(y_0)}E_{\overline
{Q}{}^{\eta}}\biggl[e^{-\varphi(Y_T)} \mathcal{E}\biggl(\int_0^{\cdot
}\varphi_y(Y_t)\rho\dd\overline{W}{}^{\eta}_t\biggr)_T\biggr].
\end{eqnarray*}
Applying once more a Girsanov transformation to eliminate the It\^{o}
exponential, we obtain
%
\begin{equation}
\label{eq:gl110}
E_{Q^{\eta}}[(X_T^{\pi^*})^\lambda]\geq x_0^\lambda e^{\widetilde
{\Lambda}(\lambda) T+\varphi(y_0)}E_{\widehat{Q}^{\eta}}\bigl[e^{-\varphi(Y_T)}\bigr],
\end{equation}
where the expectation is taken with respect to the probability measure
$\widehat{Q}^{\eta}$ on $(\Omega,\mathcal{F}_T)$ defined by
\[
\frac{d\widehat{Q}^{\eta}}{d\overline{Q}{}^{\eta}}\bigg|_{\mathcal
{F}_T}:=\mathcal{E}\biggl(\int_0^{\cdot}\varphi_y(Y_t)\rho\dd\overline
{W}{}^{\eta}_t\biggr)_T.
\]
In particular, (\ref{eq:gl110}) means that
%
\begin{equation}
\label{eq:gl416}\quad
\mathop{\underline{\lim}}_{T\uparrow\infty} \frac{1}{T}\ln
\inf_{Q^{\eta}\in\mathcal{Q}} E_{Q^{\eta}}[(X_T^{\pi^*})^\lambda]\geq
\widetilde{\Lambda}(\lambda)+\mathop{\underline{\lim}}_{T\uparrow
\infty} \frac{1}{T}\ln\inf_{Q^{\eta}\in\mathcal{Q}} E_{\widehat
{Q}^{\eta}}\bigl[e^{-\varphi(Y_T)}\bigr].
\end{equation}
Since $\vert\widetilde{h}(\pi^*(y),\eta,y)\vert^2\leq K_4(1+y^2) $,
this second change of measure can be justified again by Lemma \ref
{th:liptserersatz} combined with \cite{liptsershiryaev}, Example 3 of
Section~6.2. By Girsanov's theorem, the dynamics of $Y$ under the
new probability measure $\widehat{Q}^{\eta}$ is given by
\[
dY_t=\bigl(\widetilde{h}(\pi^*_t,\eta_t,Y_t)+\Vert\rho\Vert^2\varphi
_y(Y_t)\bigr)\dd t+\rho\dd\widehat{W}^{\eta}_t, \qquad t\leq T,
\]
where $\widehat{W}_t^{\eta}$ is a two-dimensional $\widehat{Q}^{\eta
}$-Wiener process. Moreover, inserting the definition (\ref{eq:gl24})
of $\pi^*(y)$, a straightforward computation yields the identity
\[
\widetilde{h}(\pi^*(y),\eta,y)+\Vert\rho\Vert^2\varphi_y(y)=\widetilde
{\kappa}(\eta,y).
\]
Here the function $\widetilde{\kappa}$ introduced in Assumption
\ref{th:assump4} satisfies the inequality $y\widetilde{\kappa}(\eta
$, $y)\leq-C_4y^2+C_5$ for all $\eta\in\Gamma$
with appropriate constants $C_4,C_5>0$. Thus, by
Lemma \ref{th:boundedmoments}, the quadratic moments $E_{\widehat
{Q}^{\eta}}[Y_T^2]$ are bounded above uniformly with respect to all
processes $\eta\in\mathcal{C}$ and $T\geq0$, that is,
\[
\sup_{T\geq0}\sup_{\eta\in\mathcal{C}}E_{\widehat{Q}^{\eta
}}[Y_T^2]\leq K_5(1+y_0^2).
\]
Note now that $\vert\varphi(y)\vert\leq K_6(1+y^2)$ for some constant
$K_6>0$, since the first derivative $\varphi_y$ grows at most linearly
[cf. Assumption \ref{th:assump3}(a)]. Using Jensen's
inequality, we obtain the lower bound
\begin{eqnarray*}
\ln\inf_{Q^{\eta}\in\mathcal{Q}}E_{\widehat{Q}^{\eta}}\bigl[e^{-\varphi
(Y_T)}\bigr]&\geq&\inf_{\eta\in\mathcal{C}}E_{\widehat{Q}^{\eta}}[-\varphi
(Y_T)]\geq-K_6\Bigl(1+\sup_{\eta\in\mathcal{C}}E_{\widehat{Q}^{\eta}}[
Y_T^2]\Bigr)\\
&\geq&-K_6\bigl(1+K_5(1+y_0^2) \bigr)
\end{eqnarray*}
for any finite horizon $T$. Thus the last term in (\ref{eq:gl416})
nonnegative, and so the desired estimate (\ref{eq:gl111}) follows from
(\ref{eq:gl416}).

(2) It remains to verify that the solution $(\widetilde{\Lambda}(\lambda
),\varphi)$ to our EBE (\ref{eq:gl104}) also satisfies (\ref{eq:gl106})
and vice versa. In other words, the EBE (\ref{eq:gl106}) is an
alternative version of the original equation (\ref{eq:gl104}). For this
purpose, we use the minimizing functions $\eta^*$ and $\nu^*$ defined
in Theorem \ref{th:verification1} and write $\eta^*$, $\nu^*$ and
$\pi^*$ instead of $\eta^*(y)$, $\nu^*(y)$ and $\pi^*(y)$ to simplify
the notation. Then an easy but tedious computation yields the identity
\[
\widetilde{\Lambda}(\lambda)=\tfrac{1}{2}\Vert\rho\Vert^2[\varphi
_{yy}(y)+\varphi^2_y(y)]+\widetilde{l}(\pi^*,\eta^*,y)+\varphi
_y(y)\widetilde{h}(\pi^*,\eta^*,y).
\]
Thus the pair $(\widetilde{\Lambda}(\lambda),\varphi)$ also solves the
EBE (\ref{eq:gl106}) if and only if for all $\eta\in\Gamma$
%
\begin{equation}
\label{eq:gl396}\qquad
0\leq\widetilde{l}(\pi^*,\eta,y)+\varphi_y(y)\widetilde{h}(\pi^*,\eta
,y)-[\widetilde{l}(\pi^*,\eta^*,y)+\varphi_y(y)\widetilde{h}(\pi^*,\eta^*,y)].
\end{equation}
Inserting formula (\ref{eq:gl24}) for $\pi^*$, this inequality takes
the explicit form
%
\begin{eqnarray}
\label{eq:gl107}
0&\leq&\frac{\lambda}{1-\lambda}[(\eta^{11}-\eta^{11,*})y+(\eta
^{21}-\eta^{21,*})][\theta(y)+\eta^{11,*}y+\eta^{21,*}]
\nonumber\\
&&{}+\frac{1}{1-\lambda}\rho_1\varphi_y(y)[(\eta^{11}-\eta
^{11,*})y+(\eta^{21}-\eta^{21,*})]\\
&&{}+\rho_2\varphi_y(y)[(\eta^{21}-\eta^{21,*})y+(\eta^{22}-\eta
^{22,*})]\nonumber
\end{eqnarray}
for all $\eta\in\Gamma$. To derive (\ref{eq:gl107}), we fix $\eta\in
\Gamma$ and define the convex combination $\widetilde{\eta}_\alpha:=\eta
^*+\alpha(\eta-\eta^*)$, $\alpha\in(0,1)$. Then $\widetilde{\eta}_\alpha
$ belongs to $\Gamma$, due to convexity of this set. Moreover, using
the minimizers $\eta^*,\nu^*$ and the specific choice $\nu^*_\alpha
(y):=\nu^*(\widetilde{\eta}_\alpha,y)=-\widetilde{\eta}^{12}_\alpha
y-\widetilde{\eta}_\alpha^{22}-\rho_2(y)\varphi_y(y)$, we easily derive
the inequality
\begin{eqnarray*}
0 &\leq&( 1 - \lambda) l( \widetilde{\eta}_\alpha, \nu
^*_\alpha( y ) , y ) + \varphi_y ( y ) h ( \widetilde{\eta}_\alpha, \nu
^*_\alpha( y ) , y )\\
&&{} - [ (1 - \lambda) l( \eta^* , \nu
^* , y ) + \varphi_y ( y ) h( \eta
^* , \nu^* , y ) ]
\\
&=&
\alpha[\mbox{terms in (\ref{eq:gl107})}]+\frac{1}{2}\frac{\lambda}{1-\lambda
}\alpha^2[(\eta^{11}-\eta^{11,*})y+(\eta^{21}-\eta^{21,*})]^2.
\end{eqnarray*}
Dividing finally by $\alpha$ and letting afterwards $\alpha$ tend to
zero yields the desired estimate (\ref{eq:gl107}) and equivalently (\ref
{eq:gl396}). Thus we have shown that the solution $(\widetilde{\Lambda
}(\lambda),\varphi)$ to the EBE (\ref{eq:gl104}) also satisfies (\ref
{eq:gl106}). This completes the proof.
\end{pf}
\begin{remark}
\label{th:differentialgame} The duality approach used above requires
two verification theorems. The first one characterizes the growth rate
of $U_T(x_0)$ in terms of the EBE~(\ref{eq:gl104}), and the second one
identifies an optimal long-term investment strategy and the associated
optimal growth rate $\Lambda(\lambda)$. In this remark, we discuss
\textit{heuristically} a more direct approach to (\ref{eq:gl27}) via
stochastic differential game techniques (see, e.g.,
\cite{flemingsouganidis} for an introduction). To this end, note that
(if $X_t^\pi>0$ for all $t$)
\[
E_{Q^\eta}[(X_T^\pi)^\lambda]=x_0^\lambda E_{R^{\eta,\pi}}\bigl[e^{\int
_0^T\widetilde{l}(\pi_t,\eta_t,Y_t)\dd t}\bigr],
\]
where $\widetilde{l}$ is defined in (\ref{eq:gl108}), the measure
$R^{\pi,\eta}$ is introduced in (\ref{eq:measuredifferential}), and the
dynamics of $Y$ under $R^{\eta,\pi}$ is specified in (\ref{eq:gl338}).
This suggests that
\[
U_T(x_0)=x_0^\lambda v^{u}(y,T):=x_0^\lambda\sup_{\pi\in\mathcal
{A}_T}\inf_{\eta\in\mathcal{C}} E_{R^{\eta,\pi}}\bigl[e^{\int_0^T\widetilde
{l}(\pi_t,\eta_t,Y_t)\dd t}\bigr],\qquad Y_0=y,
\]
where $v^{u}$ can be seen as the upper value function of a stochastic
differential game with maximizing ``player'' $\pi
$ and minimizing ``player'' $\eta$. The function~$v^{u}$ should be determined by the HJB-Isaacs equation
\[
v^{u}_t=\frac{1}{2}\Vert\rho\Vert^2 v^{u}_{yy}+\sup_{\pi\in\mathbb
{R}}\inf_{\eta\in\Gamma}\{\widetilde{l}(\pi,\eta,\cdot)v^{u} +\widetilde
{h}(\pi,\eta,\cdot)v^{u}_y\},\qquad v^{u}(\cdot,0)\equiv1.
\]
Using the heuristic transform $\ln v^{u}(y,T)\approx\ln U_T(x_0)\approx
\Lambda(\lambda)T+\varphi(y)$, this translates into the following EBE
of Isaacs type:
\[
\Lambda(\lambda)=\frac{1}{2}\Vert\rho\Vert^2[\varphi_{yy}+\varphi
^2_y]+\sup_{\pi\in\mathbb{R}}\inf_{\eta\in\Gamma}\{\widetilde{l}(\pi
,\eta,\cdot)+\varphi_y\widetilde{h}(\pi,\eta,\cdot)\}.
\]
If this equation has a solution $(\Lambda(\lambda),\varphi)$, then it
is easy to show that $\sup$ and $\inf$ can be interchanged and that the
saddle point is attained by $\pi^*(y)$ in (\ref{eq:gl24}) and $\eta
^*(y)$ defined in~(\ref{eq:gl410}); that is, the EBE of Isaacs type is
actually a~version of (\ref{eq:gl104}). We conjecture that the
alternative approach via differential games is also feasible. However,
the detailed derivation would be a~lenghty and technical exercise that
is beyond the scope of the present paper.\looseness=1
\end{remark}
%
\section{Existence of a solution to the ergodic Bellman equation}
\label{sec:existenceergodicBellman}
Our results rely on the existence of a specific solution $(\widetilde
{\Lambda}(\lambda),\varphi)\in\mathbb{R}_+\times C^2(\mathbb{R})$ to
the EBE (\ref{eq:gl104}). More generally, an EBE is given by
%
\begin{equation}
\label{eq:gl536}
\widetilde{\Lambda}=D\varphi(x)+H(x,\nabla\varphi)+q(x),\qquad x\in
\mathbb{R}^d,
\end{equation}
where $q$ maps from $\mathbb{R}^d$ to $\mathbb{R}$, $D$ is a second
order differential operator, and where $H$ is a real-valued nonlinear
function of the gradient $\nabla\varphi$, called the Hamiltonian. A
solution to (\ref{eq:gl536}) is a pair $(\widetilde{\Lambda},\varphi)$
of a constant $\widetilde{\Lambda}$ and a function $\varphi\dvtx\mathbb
{R}^d\rightarrow\mathbb{R}$. Such equations have been analyzed by
various authors (see, e.g., \cite{flemingmceneaney,kaisesheu,nagai}
for a discussion related to risk-sensitive
control problems, or \cite{bensoussanfrehse,bensoussannagai}).
Unfortunately their existence results do not, in general, apply to our
EBE~(\ref{eq:gl104}). The main difficulty relies on three facts: we
consider a model with nonlinear coefficients $r$, $g$ and $m$
appearing in the functions $l$ and $h$; the cost function $l$ may grow
quadratically in~$y$; (\ref{eq:gl104}) exhibits a nonlinearity with
respect to the first derivative $\varphi_y$. If the discussion is
limited to linear coefficients, then a quadratic ansatz may yield an
explicit solution to (\ref{eq:gl104}) (see, e.g.,
\cite{flemingsheu,pham}, and also Section \ref{sec:geometricOUutility} for a
case study).

Let us now turn to the existence problem for nonlinear coefficients.
The EBE (\ref{eq:gl104}) can be rewritten in the condensed form
%
\begin{equation}
\label{eq:gl41}\quad
\widetilde{\Lambda}(\lambda)=\frac{1}{2}\Vert\rho\Vert^2\varphi
_{yy}(y)+\frac{1}{2}(\widehat{\rho}\varphi_y(y))^2+\inf_{\eta\in\Gamma
}\{n(\eta,y)+\varphi_y(y)m(\eta,y)\},
\end{equation}
where we use the notation $\widehat{\rho}:=\sqrt{\frac{1}{1-\lambda
}\rho_1^2+\rho_2^2}$,
\begin{eqnarray*}
n(\eta,y)&:=&\frac{1}{2}\frac{\lambda}{1-\lambda}[\theta(y)+\eta
^{11}y+\eta^{21}]^2+\lambda r(y),\\
m(\eta,y)&:=&g(y)+\frac{1}{1-\lambda
}\rho_1\bigl(\lambda\theta(y)+\eta^{11}y+\eta^{21}\bigr)+\rho_2(\eta^{12}y+\eta^{22}).
\end{eqnarray*}
The following existence result is deduced from Fleming and McEneaney
\cite{flemingmceneaney}. Their construction of a solution involves a
parameterized family of finite time horizon \textit{stochastic
differential games} (see, e.g., Fleming and Souganidis \cite
{flemingsouganidis}). The associated value function is characterized in
terms of a parabolic PDE, called \textit{Isaacs' equation}, and the
existence of a solution ($\widetilde{\Lambda}(\lambda),\varphi$)
follows by taking appropriate limits of the Isaacs' PDE when both
``time'' tends to infinity and the underlying
parameter converges to zero.
\begin{lemma}
In addition to Assumption \ref{th:assump1} let us assume that
$\theta$ is bounded, that $\Gamma\subset\{(0,0)\}\times\mathbb{R}^2$
and that
%
\begin{equation}
\label{eq:gl173}
\exists K>0\dvtx g_y(y)+\frac{\lambda}{1-\lambda
}\rho_1\theta_y(y)\leq-K \qquad\mbox{for all } y\in\mathbb{R}.
\end{equation}
Then there exist a pair $ \widetilde{\Lambda}(\lambda)\in\mathbb{R}_+$,
$\varphi\in C^2(\mathbb{R})$ that solves the EBE (\ref{eq:gl104}).
Moreover, we have $\vert\varphi_y\vert\leq\max_{\eta\in\Gamma}\Vert
n_y(\eta,\cdot)\Vert_\infty/K$, and so this solution also satisfies the
regularity Assumptions \ref{th:assump3} and \ref{th:assump4}.
\end{lemma}
\begin{pf} Our assumptions ensure boundedness of $n\geq0$, $\eta_y$
and $m_y$ on $\Gamma\times\mathbb{R}$. Moreover, the mean value
theorem combined with (\ref{eq:gl173}) gives
\[
(x-y)\bigl(m(\eta,x)-m(\eta,y)\bigr)\leq-K\vert x-y \vert^2 \qquad\mbox{for
all } x,y\in\mathbb{R},\eta\in\Gamma.
\]
The functions $n$, $m$ thus satisfy condition (7.2) in Fleming and McEneaney~\cite{flemingmceneaney}, and applying \cite{flemingmceneaney},
Theorem 7.1, for $\gamma:=(\sqrt{2}\widehat{\rho})^{-1}$ and
$\varepsilon
:=\Vert\rho\Vert^2/\widehat{\rho}^2$ the desired existence result
follows.
\end{pf}
%
\section{Explicit results}
\label{sec:explicitresults}
\subsection{Black--Scholes model with uncertain drift}
For constant coefficients $r(y)\equiv r$ and $m(y)\equiv m$, the
reference model $Q_0$ in Section \ref{sec:utilitymax} becomes the
Black--Scholes model with price dynamics
\[
dS^0_t=S^0_t r\dd t,\qquad dS^1_t=S^1_t(m\dd t+\sigma\dd W^{1}_t).
\]
In particular, the market price of risk function $\theta(y)=\frac
{m-r}{\sigma}$ is constant. Taking the specific set $\Gamma=\{(0,0)\}
\times[a,b]\times\{0\}$, $a\leq0\leq b$, each measure $Q^{\eta}\in
\mathcal{Q}$ corresponds to a drift perturbation of the following type:
\[
dS^1_t=S^1_t\bigl([m+\sigma\eta^{21}_t]\dd t+\sigma\dd W^{1,\eta}_t\bigr).
\]
In this example the factor process $Y$ plays no role. In particular,
the maximal expected utility for a finite horizon does not depend on
the initial state of the factor process. Hence the function $\varphi$
appearing in the heuristic separation of time and space variables (\ref
{eq:gl249}) is constant, and its derivatives $\varphi_y$, $\varphi
_{yy}$ vanish. The EBE (\ref{eq:gl104}) thus reduces to
\[
\widetilde{\Lambda}(\lambda)=\inf_{\nu\in\mathbb{R}}\inf_{\eta\in\Gamma
}\biggl\{\frac{1}{2}\frac{\lambda}{1-\lambda}[(\theta+\eta^{21})^2+\nu
^2]+\lambda r\biggr\}=\frac{1}{2}\frac{\lambda}{1-\lambda}\inf_{\eta\in
\Gamma}\{(\theta+\eta^{21})^2\}+\lambda r.
\]
The number $\widetilde{\Lambda}(\lambda)$ can be expressed in terms of
the element $\eta^{21,*}\in[a,b]$ which minimizes the absolute value
$\vert\theta+\eta^{21}\vert$ among all $\eta^{21}\in[a,b]$. Defining
the constant controls
$\eta_t^*:=(0,0,\eta^{21,*},0)$ and $\nu_t^*:=0$, $t\geq0$,
the verification theorems can be transferred to our present example in
a simplified form which does not not require any additional conditions
as in Assumptions \ref{th:assump3} and \ref{th:assump4}. As a
result we get the following description of the aymptotics of robust
expected power utility:
\begin{itemize}
\item[$\bullet$] The maximal robust utility $U_T(x_0)$ grows
exponentially with rate
\[
\widetilde{\Lambda}(\lambda)=\frac{1}{2}\frac{\lambda}{1-\lambda
}(\theta+\eta^{21,*})^2+\lambda r>0.
\]
\item[$\bullet$] $\Lambda(\lambda)=\sup_{\pi\in\mathcal{A}}
\overline{\lim}_{T\uparrow\infty} \frac{1}{T}\ln
\inf_{Q^{\eta}\in\mathcal{Q}} E_{Q^{\eta}}[(X_T^{\pi})^\lambda
]=\widetilde{\Lambda}(\lambda)$.
\item[$\bullet$] The optimal long-term strategy takes the form
\[
\pi^*_t:=\frac{1}{1-\lambda}\frac{1}{\sigma}(\theta+\eta
^{21,*}),\qquad t\geq0.\vadjust{\goodbreak}
\]
\item[$\bullet$]
The asymptotic worst-case model $Q^{\eta^*}$ is given by the constant
control $\eta_t^*=(0,0,\eta^{21,*},0)$, and it does not depend on the
parameter $\lambda$.
\end{itemize}
\begin{remark} Using methods from robust statistics, Schied \cite
{schied} shows that the measure $Q^{\eta^*}$ is actually \textit{least
favorable} in the following sense: for any finite maturity, the robust
utility maximization problem (\ref{eq:gl295}) is equivalent to the
classical problem for $Q^{\eta^*}$, regardless of the choice of the
underlying utility function $u$.
\end{remark}
%
\subsection{Geometric Ornstein--Uhlenbeck model with uncertain mean reversion}
\label{sec:geometricOUutility}
As our second case study, we consider the case where the economic
factor $Y$ is of OU type, and where there interest rate $r$ is
constant. In our reference model $Q_0$, the factor $Y$ is assumed to be
a classical OU process with constant rate of mean reversion $\eta_0>0$
and volatility $\sigma>0$, that is,
%
\begin{equation}
\label{eq:gl75}
dY_t=-\eta_0 Y_t\dd t+\sigma\dd W^1_t,\qquad Y_0=y_0.
\end{equation}
We assume that $S^1_t:=\exp(Y_t+\alpha t)$, $\alpha\in\mathbb{R}$,
describes the price process of the risky asset. By It\^{o}'s
formula, the dynamics of $S^1$ is governed by the SDE
\[
dS^1_t=S^1_t\bigl(\alpha\dd t+dY_t+\tfrac{1}{2}\dd\langle Y\rangle
_t\bigr)=S^1_t\bigl(\bigl(-\eta_0 Y_t+\tfrac{1}{2}\sigma^2+\alpha\bigr)\dd t+\sigma\dd W^1_t\bigr).
\]
Hence this example corresponds to the general model of Section \ref
{sec:utilitymax} for the choice $g(y)=-\eta_0y$, $\rho_1=\sigma$, $\rho
_2=0$, $m(y)=-\eta_0y+\frac{1}{2}\sigma^2+\alpha$, and for the affine
market price of risk function
\[
\theta(y)=\frac{1}{\sigma}\biggl(-\eta_0y+\frac{1}{2}\sigma^2+\alpha-r\biggr).
\]
Let us suppose that the investor is uncertain about the ``true''
future rate of mean reversion: instead of a constant
rate we admit any rate process that is progressively measurable and
that takes its values in some interval $[a,b]$, $0<a\leq\eta_0\leq
b<\infty$. This uncertainty about the true rate of mean reversion can
be embedded into our general model by choosing the set
\[
\Gamma=\biggl[\frac{\eta_0-b}{\sigma},\frac{\eta_0-a}{\sigma}\biggr]\times\{
(0,0,0)\}.
\]
Indeed, let $Q^{\eta}\in\mathcal{Q}$ denote the probabilistic model
generated by a $\Gamma$-valued, progressively measurable process $\eta
=(\eta_t)_{t\geq0}$; cf. (\ref{eq:gl381}). In view of (\ref
{eq:gl267}), the factor process $Y$ then evolves under $Q^{\eta}$
according to
\[
dY_t=-(\eta_0-\sigma\eta^{11}_t)Y_t\dd t+\sigma\dd W_t^{1,\eta},
\]
and the resulting mean reversion process $(\eta_0-\sigma\eta
^{11}_t)_{t\geq0}$ takes values in $[a,b]$.

To prepare the analysis of the asymptotic robust utility maximization
problem~(\ref{eq:gl27}), we first solve its nonrobust version
%
\begin{equation}
\label{eq:gl312}
\mbox{maximize } \mathop{\overline{\lim}}_{T\uparrow\infty
} \frac{1}{T}\ln E_{Q_0}[(X_T^{\pi})^\lambda] \mbox{ among
all } \pi\in\mathcal{A}
\end{equation}
for the specific model $Q_0$. This problem has been studied, amongst
others, by Fleming and Sheu \cite{flemingsheu} and Pham \cite{pham}. By
the following proposition we recover their results as a special case of
our general robust duality approach. To indicate the nonrobust case,
we denote the optimal growth rate for (\ref{eq:gl312}) by $\Lambda
_{Q_0}(\lambda)$. Note that $\mathcal{Q}=\{Q_0\}$ if we take the
one-point set $\Gamma=\{(0,0,0,0)\}$. Thus our general EBE (\ref
{eq:gl104}) takes the simplified form
%
\begin{eqnarray}
\label{eq:gl502}
\widetilde{\Lambda}_{Q_0}(\lambda)
&=&\frac{1}{2}\sigma^2\biggl[\varphi_{yy}(y)+\frac{1}{1-\lambda}\varphi
_y^2(y)\biggr]+\lambda r\nonumber\\
&&{}+\frac{1}{2}
\frac{\lambda}{1-\lambda}\biggl(\frac{-\eta_0y+({1/2})\sigma^2+\alpha
-r}{\sigma}\biggr)^2\\
&&{}+\varphi_y(y)\biggl[-\frac{1}{1-\lambda}\eta_0y+\frac{\lambda}{1-\lambda
}\biggl(\frac{1}{2}\sigma^2+\alpha-r\biggr)\biggr],\nonumber
\end{eqnarray}
where the infimum among all $\nu\in\mathbb{R}$ is attained for $\nu
^*(y)\equiv0$.
\begin{proposition}
\label{th:OUlongterm}
The EBE (\ref{eq:gl502}) has the solution
%
\begin{subequation}
\label{eq:gl382}
\begin{eqnarray}
\widetilde{\Lambda}_{Q_0}(\lambda)&=&\frac{1}{2}\bigl(1-\sqrt{1-\lambda
}\bigr)\eta_0+\lambda\biggl(r+\frac{1}{2\sigma^2}\biggl(\frac{1}{2}\sigma^2+\alpha
-r\biggr)^2\biggr),\\
\varphi(y)&=&\frac{1}{2}\bigl(1-\sqrt{1-\lambda}\bigr)\frac{\eta
_0}{\sigma^2}y^2-\frac{\lambda}{\sigma^2}\biggl(\frac{1}{2}\sigma^2+\alpha-r\biggr)y,
\end{eqnarray}
\end{subequation}
which satisfies our regularity Assumptions \ref{th:assump3} and
\ref{th:assump4}. Thus it holds that
\[
\widetilde{\Lambda}_{Q_0}(\lambda)=\lim_{T\uparrow\infty}\frac
{1}{T}\ln U_T^{Q_0}(x_0)=\Lambda_{Q_0}(\lambda)=\sup_{\pi\in\mathcal
{A}}\mathop{\overline{\lim}}_{T\uparrow\infty} \frac
{1}{T}\ln E_{Q_0}[(X_T^{\pi})^\lambda].
\]
Moreover, an optimal feedback strategy $\pi^*_t=\pi^*(Y_t)$, $t\geq0$,
for our investment problem (\ref{eq:gl312}) is given by the affine function
%
\begin{equation}
\label{eq:gl503}
\pi^*(y)=-\frac{1}{\sqrt{1-\lambda}}\frac{\eta_0}{\sigma^2}y+\frac
{1}{\sigma^2}\biggl(\frac{1}{2}\sigma^2+\alpha-r\biggr).
\end{equation}
\end{proposition}
\begin{pf} Following \cite{flemingsheu} and \cite{pham} we are looking
for a quadratic solution $
\varphi(y)=\frac{1}{2}Ay^2+By$. Inserting the derivatives in (\ref
{eq:gl502}) and comparing the coefficients of the terms in $y^2$, in
$y$, and the constants yields that the EBE (\ref{eq:gl502}) holds for
every triple $(A,B,\widetilde{\Lambda}_{Q_0}(\lambda))$ satisfying the
system of equations
\begin{eqnarray*}
0&=&\frac{1}{2}\sigma^2 A^2-\eta_0A+\frac{\lambda}{2\sigma^2} \eta
_0^2,\\
0&=&\sigma^2AB+\lambda\biggl(\frac{1}{2}\sigma^2+\alpha-r\biggr)A-B\eta_0-\frac
{\lambda}{\sigma^2}\biggl(\frac{1}{2}\sigma^2+\alpha-r\biggr)\eta_0, \\
\widetilde{\Lambda}_{Q_0}(\lambda)&=&\frac{1}{2}\sigma^2 \biggl(A+\frac
{1}{1-\lambda}B^2\biggr)+\frac{\lambda}{1-\lambda}\biggl(\frac{1}{2}\sigma
^2+\alpha-r\biggr)B+\lambda r\\
&&{}+\frac{1}{2}\frac{\lambda}{1-\lambda}\biggl(\frac
{({1/2})\sigma^2+\alpha-r}{\sigma}\biggr)^2.
\end{eqnarray*}
The quadratic equation for $A$ has the solutions $A_{\pm}=(1\pm\sqrt
{1-\lambda})\frac{\eta_0}{\sigma^2}$. We choose $A=A_-$, and we shall
explain in Remark \ref{th:irrelevant} why the other solution is
irrelevant. A straightforward calculation gives $B=-\frac{\lambda
}{\sigma^2}(\frac{1}{2}\sigma^2+\alpha-r)$ and finally the expressions
for $\widetilde{\Lambda}_{Q_0}(\lambda)$ and $\varphi$ in (\ref
{eq:gl382}). The parabola $\varphi$ is bounded below, $\varphi_y$ grows
linearly and the functions $\kappa$, $\widetilde{\kappa}$ defined in
Assumption \ref{th:assump3}(b) and \ref{th:assump4} satisfy the
regularity condition
\[
y\kappa(0,y)=y\widetilde{\kappa}(0,y)=-\frac{1}{\sqrt{1-\lambda}}\eta_0y^2.
\]
Applying Theorems \ref{th:verification1} and \ref
{th:identification} completes the proof.
\end{pf}
\begin{remark}
\label{th:irrelevant}
Using the other root $A_+$ yields $\varphi(y)=\frac{1}{2}A_+ y^2+By$ and
\[
\widetilde{\Lambda}_{Q_0}(\lambda)=\frac{1}{2}\bigl(1+\sqrt{1-\lambda}\bigr)\eta
_0+\lambda\biggl(r+\frac{1}{2\sigma^2}\biggl(\frac{1}{2}\sigma^2+\alpha-r\biggr)^2\biggr).
\]
In particular, this example illustrates that the solution to an EBE is
not necessarily unique. On the other hand, the ``ergodicity''
Assumption \ref{th:assump3}(b) selects the good
candidate. Indeed, the proof of Theorem \ref{th:verification1}
requires that $\lim_{T\uparrow\infty}\frac{1}{T}E_{\widehat{R}^{\eta
}}[Y_T^2]=0$, where $Y$ follows the SDE (\ref{eq:gl157}). Given the
geometric OU model and the solutions $\varphi(y)=\frac{1}{2}A_\pm
y^2+By$ this SDE takes the form
\[
dY_t=\pm\frac{\eta_0}{\sqrt{1-\lambda}}Y_t\dd t+ \sigma
d\widehat{W}^{1,\eta}_t,\qquad Y_0=y_0.
\]
The factor process $Y$ is an ``explosive''
Gaussian process for the root $A_+$ in the sense that $\lim
_{T\uparrow\infty}\frac{1}{T}E_{\widehat{R}^{\eta}}[Y_T^2]=\infty$.
Thus the arguments used in the proof of Theorem \ref
{th:verification1} fail and so the solution associated with $A_+$ is
irrelevant. Conversely, taking $A_-$, Theorem \ref
{th:boundedmoments} applies to the ergodic process $Y$.
\end{remark}

As a complement to Proposition \ref{th:OUlongterm} we look at the
maximal robust utility~$U^{Q_0}_T(x_0)$ attainable at time $T$ and the
asymptotics of the optimal investment strategy $\pi^{*,T}$ as $T\uparrow
\infty$. The following proposition extends Propositions~5.6 and
5.7 in F\"{o}llmer and Schachermayer \cite{foellmerschachermayer} by
including an interest rate $r>0$ and the additional drift component
$\alpha$ for the price process $S^1$.
\begin{proposition}
\label{th:OUutilityfinite}
For any initial condition $Y_0=y_0$, the maximal robust expected
utility $U^{Q_0}_T(x_0)$ takes the form
%
\begin{equation}
\label{eq:gl77}
\qquad U^{Q_0}_T(x_0)=\frac{1}{\lambda}x_0^\lambda
\bigl[(A_T^-)^{-1/2}e^{B_T(y_0)+({\lambda}/({1-\lambda
}))rT+(A_T^-)^{-1}C_T(y_0)}\bigr]^{1-\lambda},
\end{equation}
where we use the notation
\begin{eqnarray*}
A_T^\pm&:=&1-\tfrac{1}{2}\bigl(1-(1-\lambda)^{-1/2}\bigr)\bigl(1\pm\exp\bigl(-2\eta
_0(1-\lambda)^{-1/2}T\bigr)\bigr),\\
B_T(y)&:=&-\frac{\eta_0}{2\sigma^2}[(1-\lambda)^{-1/2}-(1-\lambda
)^{-1}]y^2+\frac{1}{\sigma^2}\frac{\lambda}{\lambda-1}\biggl(\frac
{1}{2}\sigma^2+\alpha-r\biggr)y\\
&&{}-\frac{1}{2}\biggl[\eta
_0\bigl((1-\lambda)^{-1/2}-(1-\lambda)^{-1}\bigr)+\frac{1}{\sigma^2}\frac
{\lambda}{\lambda-1}\biggl(\frac{1}{2}\sigma^2+\alpha-r\biggr)^2\biggr]T,\\
C_T(y)&:=&\frac{\eta_0}{2\sigma^2}\bigl((1-\lambda)^{-1/2}-(1-\lambda
)^{-1}\bigr)\exp\bigl(-2\eta_0(1-\lambda)^{-1/2}T\bigr)y^2\\
&&{}-\frac
{1}{\sigma^2}\frac{\lambda}{\lambda-1}\biggl(\frac{1}{2}\sigma^2+\alpha
-r\biggr)\exp\bigl(-\eta_0(1-\lambda)^{-1/2}T\bigr)y\\
&&{}+\frac{1}{4\sigma
^2}\frac{\lambda^2}{(1-\lambda)^{3/2}}\biggl(\frac{1}{2}\sigma^2+\alpha
-r\biggr)^2\bigl(1-\exp\bigl(-2\eta_0(1-\lambda)^{-1/2}T\bigr)\bigr).
\end{eqnarray*}
The optimal proportion $\pi_t^{*,T}$ is an affine function of the
current state $Y_t$ of the factor process given by $\pi
_t^{*,T}=a[T-t]Y_t+b[T-t]$, where
\begin{eqnarray*}
a[T-t]&:=&-\frac{\eta_0}{\sigma^2}(1-\lambda
)^{-1/2}A^+_{T-t}(A^-_{T-t})^{-1},\\
b[T-t]&:=&\frac{1}{\sigma^2}\biggl(\frac{1}{2}\sigma^2+\alpha
-r\biggr)\biggl[1+(A_{T-t}^-)^{-1}\frac{\lambda}{1-\lambda}e^{-\eta_0(1-\lambda
)^{-1/2}(T-t)}\biggr].
\end{eqnarray*}
\end{proposition}
\begin{pf}
Detailed computations can be found in \cite{knispel}, Chapter 4.
\end{pf}

Since $A^\pm_T$ and $C_T(y_0)$ converge to a finite limit as $T\uparrow
\infty$, we thus obtain
\[
\lim_{T\uparrow\infty}\frac{1}{T}\ln U^{Q_0}_T(x_0)=(1-\lambda)\lim
_{T\uparrow\infty}\frac{1}{T}\biggl(B_T(y_0)+\frac{\lambda}{1-\lambda
}rT\biggr)=\widetilde{\Lambda}_{Q_0}(\lambda),
\]
in accordance with Proposition \ref{th:OUlongterm}. Moreover, we have
\[
\lim_{T\uparrow\infty}a[T-t]=-\frac{\eta_0}{\sigma^2\sqrt{1-\lambda
}} \quad\mbox{and}\quad \lim_{T\uparrow\infty
}b[T-t]=\frac{1}{\sigma^2}\biggl(\frac{1}{2}\sigma^2+\alpha-r\biggr),
\]
due to $\lim_{T\uparrow\infty}A^{\pm}_T=\frac{1}{2}(1+(1-\lambda)^{-1/2})$.
Thus the asymptotic form of the optimal strategy $\pi^{*,T}$ as
$T\uparrow\infty$ is given by
%
\begin{equation}
\lim_{T\uparrow\infty}\pi_t^{*,T}=-\frac{1}{\sigma^2\sqrt{1-\lambda
}}\eta_0Y_t+\frac{1}{\sigma^2}\biggl(\frac{1}{2}\sigma^2+\alpha-r\biggr),
\end{equation}
and so it coincides with the optimal long-term strategy $\pi^*$ in (\ref
{eq:gl503}). On the other hand, Fleming and Sheu \cite{flemingsheu}
observed that $\lim_{T\uparrow\infty}\pi_t^{*,T}$ does not provide an
optimal long-term strategy for power utility with parameter $\lambda\leq-3$.

Let us now analyze the \textit{robust} case. Since $\Lambda_{Q_0}(\lambda
)$ is increasing in $\eta_0$, it is natural to expect that the
asymptotic worst-case measure $Q^{\eta^*}$ corresponds to the
probabilistic model, under which $Y$ has the minimal rate of mean
reversion $a$. The following proposition confirms this conjecture.
\begin{proposition}
\label{th:robustpowerutilityOU}
For the geometric OU model with uncertain rate of mean reversion, the
optimal growth rate of robust power utility is given by
\[
\Lambda(\lambda)=\frac{1}{2}\bigl(1-\sqrt{1-\lambda}\bigr)a+\lambda\biggl(r+\frac
{1}{2\sigma^2}\biggl(\frac{1}{2}\sigma^2+\alpha-r\biggr)^2\biggr)>0,
\]
and the maximal robust utility $U_T(x_0)$ grows exponentially at this
rate. The asymptotic worst-case model $Q^{\eta^*}$ is determined by
$\eta_t^*=(\frac{\eta_0-a}{\sigma},0,0,0)$, and the optimal long term
strategy $\pi_t^*=\pi^*(Y_t)$ is specified by the affine function
\[
\pi^*(y)=-\frac{1}{\sqrt{1-\lambda}}\frac{a}{\sigma^2}y+\frac
{1}{\sigma^2}\biggl(\frac{1}{2}\sigma^2+\alpha-r\biggr).
\]
\end{proposition}
\begin{pf} Replacing $\eta_0$ in (\ref{eq:gl382}) by the minimal mean
reversion $a$ provides
\begin{eqnarray*}
\widetilde{\Lambda}(\lambda)&=&\frac{1}{2}\bigl(1-\sqrt{1-\lambda}\bigr)a+\lambda
\biggl(r+\frac{1}{2\sigma^2}\biggl(\frac{1}{2}\sigma^2+\alpha-r\biggr)^2\biggr)>0,\\
\varphi
(y)&=&\frac{1}{2}\bigl(1-\sqrt{1-\lambda}\bigr)\frac{a}{\sigma^2}y^2-\frac
{\lambda}{\sigma^2}\biggl(\frac{1}{2}\sigma^2+\alpha-r\biggr)y
\end{eqnarray*}
as a candidate for the solution to the EBE (\ref{eq:gl104}), and it is
easy to verify that~$(\widetilde{\Lambda}(\lambda),\varphi)$ is indeed
a solution. The corresponding minimizers are $\nu^*(y)\equiv0$, $\eta
^*(y)\equiv(\frac{\eta_0-a}{\sigma},0,0,0)\in\Gamma$.
It remains to verify that $(\widetilde{\Lambda}(\lambda),\varphi)$
satisfies our Assumptions \ref{th:assump3} and \ref{th:assump4}:
since $\varphi\in C^2(\mathbb{R})$ is a parabola, it is bounded below
and its first derivative $\varphi_y$ grows linearly. Moreover, the
auxiliary functions~$\kappa$,~$\widetilde{\kappa}$ appearing in our
Assumptions \ref{th:assump3} and \ref{th:assump4} satisfy for all
$\eta\in\Gamma$
\[
y\kappa(\eta,y)=\biggl[-\frac{1}{1-\lambda}(\eta_0-a-\sigma\eta^{11})-\frac
{1}{\sqrt{1-\lambda}}a\biggr]y^2\leq-\frac{1}{\sqrt{1-\lambda}}ay^2
\]
and
\[
y\widetilde{\kappa}(\eta,y)=\biggl[-(\eta_0-a)+\sigma\eta^{11}-\frac
{1}{\sqrt{1-\lambda}}a\biggr]y^2\leq-\frac{1}{\sqrt{1-\lambda}}ay^2
\]
due to $\eta^{11}\leq\frac{\eta_0-a}{\sigma}$ and $\eta^{11,*}=\frac
{\eta_0-a}{\sigma}$. We thus derive Proposition \ref
{th:robustpowerutilityOU} as a~special case of Theorems \ref
{th:verification1} and \ref{th:identification}.
\end{pf}
%
\section{Application to a robust outperformance criterion}
\label{sec:outperformance}
Utility maximization is conceptually related to specific numerical
representations of the investor's preferences. The application
requires us to know the utility function~$u$\vadjust{\goodbreak} which is by nature
subjective. For institutional managers utility maximization thus
creates severe difficulties. On the one hand, the preferences of their
customers and the corresponding numerical representations are not
really known exactly. On the other hand, the individual preferences of
the managers and of the various customers with shares in the same
investment fund will typically be different. This suggests that we
should look for an ``intersubjective'' criterion
for optimal portfolio management which is acceptable for a large
variety of investors. Such an alternative consists of evaluating the
performance of the portfolio relative to a given \textit{benchmark} such
as a~stock index. The investor aims at outperforming the benchmark with
maximal probability. If the benchmark is a contingent claim $H$ at
a~terminal time~$T$, then the outperformance problem reduces to
maximizing the probability $Q[X^\pi_T\geq H]$ of a \textit{successful
hedge}. This criterion, known as \textit{quantile hedging}, has been
developed as a substitute for investors who are not willing or not able
to raise the initial costs required by a \textit{perfect hedging} or \textit
{superhedging} strategy of $H$ (see, e.g., F\"{o}llmer and Leukert
\cite{foellmerleukert} and the references therein).

Pham \cite{pham} proposed an asymptotic benchmark criterion for
optimal long-term investment. Here the investor has in mind a level of
return $c$ and aims at maximizing the probability that the portfolio's
growth rate
\[
L_T^\pi:=\frac{1}{T}\ln X_T^\pi
\]
[or more generally $\frac{1}{T}\ln(X_T^\pi/I_T)$ for an index process
$I$] exceeds this threshold. For finite $T$ this corresponds to
quantile hedging for $H=\exp(cT)$. But what happens in the long run?
If the growth rates $L_T^\pi$ converge $Q$-a.s. as $T\uparrow\infty$
and satisfy under $Q$ a \textit{large deviations principle} with rate
function~$I^\pi$, then $Q[L_T^\pi\geq c]\approx\exp(-I^\pi(c)T)$ as
$T\uparrow\infty$; that is, the probability that~$L_T^\pi$ departs from
its limiting value decays to zero exponentially fast. Thus the long
term view amounts to minimizing the rates $I^\pi(c)$, or equivalently
to\looseness=1
%
\begin{equation}
\label{eq:gl558}
\mbox{maximizing } \mathop{\overline{\lim}}_{T\uparrow\infty
} \frac{1}{T}\ln Q[L_T^\pi\geq c] \mbox{ among all } \pi.
\end{equation}\looseness=0
An asymptotic benchmark criterion of this form may be of particular
interest for institutional managers with long-term horizon, such as
mutual fund managers. Note, however, that this ansatz does not take
into account the size of the shortfall if it does occur. From a
mathematical point of view, it leads to a \textit{large
deviations control problem}. On the other hand, standard results from
the large deviations theory (such as the G\"{a}rtner--Ellis theorem;
see, e.g., \cite{dembozeitouni}, Theorem 2.3.6) suggest that the rate
function $I^\pi$ is a~Fenchel--Legendre transform of the logarithmic
moment generating function\looseness=1
\[
\Lambda_Q(\lambda,\pi):=\mathop{\overline{\lim}}_{T\uparrow\infty
} \frac{1}{T}\ln E_Q[\exp(\lambda TL_T^\pi)]=\mathop
{\overline{\lim}}_{T\uparrow\infty} \frac{1}{T}\ln E_Q[(X_T^\pi
)^\lambda].
\]\looseness=0
In this spirit, Pham developed a duality approach to (\ref{eq:gl558}).
His Theorem~3.1, relying on large deviations arguments, but not
on the specific structure of the underlying market model, states that
%
\begin{equation}
\label{eq:gl556}
\sup_{\pi}\mathop{\overline{\lim}}_{T\uparrow\infty} \frac
{1}{T}\ln Q[L_T^\pi\geq c]=-\sup_{\lambda\in(0,\lambda')}\{\lambda
c-\Lambda_Q(\lambda)\},
\end{equation}
where $\Lambda_Q(\lambda):=\sup_{\pi}\Lambda_Q(\lambda,\pi)$ is the
optimal growth rate of expected power utility with respect to $Q$.
Applications of Pham's theorem to specific market models can be
found in \cite{pham,pham2,hataiida,hatasekine,sekine}.

However, the benchmark criterion (\ref{eq:gl558}) does not account for
model ambiguity. To overcome this limitation, it is natural to study
its robust version,
%
\begin{equation}
\label{eq:gl560}
\mbox{maximize } \mathop{\overline{\lim}}_{T\uparrow\infty
} \frac{1}{T}\ln\inf_{Q\in\mathcal{Q}} Q[L_T^\pi\geq c]
\mbox{ among all } \pi.
\end{equation}
The solution is derived in \cite{knispel}, Chapter 6, for the \textit
{robust} stochastic factor model of Section \ref{sec:utilitymax}, and
it is closely related to the asymptotics of robust utility
maximization. Under suitable regularity assumptions [e.g., $\Lambda
\in C^1((0,1))$] and $\lim_{\lambda\uparrow\lambda'}\Lambda'(\lambda
)=\infty$ for some $\lambda'\leq1$ we obtain the duality formula
%
\begin{equation}
\label{eq:gl510}
\sup_{\pi\in\mathcal{A}}\mathop{\overline{\lim}}_{T\uparrow\infty
} \frac{1}{T}\ln\inf_{Q^\eta\in\mathcal{Q}}Q^\eta[L_T^\pi\geq
c]=-\sup_{\lambda\in(0,\lambda')}\{\lambda c-\Lambda(\lambda)\}.
\end{equation}
This can be seen as a robust extension of (\ref{eq:gl556}), but here
the duality formula involves the optimal growth rates $\Lambda(\lambda
)$, $\lambda\in(0,1)$, of robust power utility. Moreover, the sequence
of investment processes $\widehat{\pi}^{c,n}$, $n\in\mathbb{N}$,
defined by
\[
\widehat{\pi}_t^{c,n}=\cases{
\pi^*_t(\lambda[c+1/n]), &\quad for $c>\Lambda'(0)$,\cr
\pi^*_t\bigl(\lambda[\Lambda'(0)+1/n]\bigr), &\quad for $c\leq
\Lambda'(0)$,}
\]
in terms of the optimal long term strategies $\pi^*(\lambda)$ for
robust power utility, and in terms of parameters $\lambda[c]\in
\arg\max_{\lambda\in(0,\lambda')}\{\lambda c-\Lambda
(\lambda)\}$ is nearly optimal for (\ref{eq:gl560}). The proof is
beyond the scope of this paper and therefore omitted.
\begin{example} For the geometric OU model with uncertain mean
reversion (see Section \ref{sec:geometricOUutility})
Proposition \ref{th:robustpowerutilityOU} shows that:
\begin{itemize}
\item[$\bullet$] $\Lambda(\lambda)=\frac{1}{2}(1-\sqrt{1-\lambda
})a+\lambda\gamma$ with $\gamma:=r+\frac{1}{2\sigma^2}(\frac
{1}{2}\sigma^2+\alpha-r)^2$,
\item[$\bullet$] $\pi_t^*(\lambda)=-\frac{1}{\sqrt{1-\lambda}}\frac
{a}{\sigma^2}Y_t+\frac{1}{\sigma^2}(\frac{1}{2}\sigma^2+\alpha-r)$,
$t\geq0$.
\end{itemize}
We thus obtain from (\ref{eq:gl510}) the optimal rate of exponential decay
\[
\sup_{\pi\in\mathcal{A}}\mathop{\overline{\lim}}_{T\uparrow\infty
} \frac{1}{T}\ln\inf_{Q\in\mathcal{Q}}Q[L_T^\pi\geq c]=\cases{
\displaystyle -\frac{({a}/{4}-c+\gamma)^2}{c-\gamma}, &\quad
for $\displaystyle c>\frac{a}{4}+\gamma$,\vspace*{3pt}\cr
0, &\quad for $\displaystyle c\leq\frac{a}{4}+\gamma$.}
\]
Since $\lambda[c]=1-(\frac{a}{4(c-\gamma)})^2$, the nearly optimal
strategies are given by
\[
\widehat{\pi}_t^{c,n}=\cases{
\displaystyle -\frac{4}{\sigma^2}\biggl(c+\frac
{1}{n}-\gamma\biggr)Y_t+\frac{1}{\sigma^2}\biggl(\frac{1}{2}\sigma^2+\alpha-r\biggr),&\quad
for $\displaystyle c>\frac{a}{4}+\gamma$,\vspace*{3pt}\cr
\displaystyle -\frac{4}{\sigma^2}\biggl(\frac{a}{4}+\frac
{1}{n}\biggr)Y_t+\frac{1}{\sigma^2}\biggl(\frac{1}{2}\sigma^2+\alpha-r\biggr), &\quad
for $\displaystyle c\leq\frac{a}{4}+\gamma$.}
\]
\end{example}
\begin{remark}
Another natural problem is to minimize the robust large deviations
probability of downside risk
%
\begin{equation}
\label{eq:gl511}
\mathop{\underline{\lim}}_{T\uparrow\infty} \frac{1}{T}\ln
\sup_{Q\in\mathcal{Q}}Q[L_T^\pi\leq c].
\end{equation}
Here the investor is interested in minimizing, in the long run, the
worst-case probability that his portfolio underperforms a savings account with
interest rate $c$. In the nonrobust case, this large deviation
criterion has been proposed by Pham \cite{pham}, but a rigorous
solution was given first by Hata, Nagai and Sheu~\cite{hatanagaisheu}
for the special case of a linear Gaussian factor model. The solution
can be derived by a duality approach which, in contrast to~(\ref
{eq:gl556}) and~(\ref{eq:gl510}), involves the optimal growth rates
$\Lambda(\lambda)$ of power utility with negative parameter~$\lambda$
(cf. Remark \ref{th:negpow}). For a detailed discussion of
problem~(\ref{eq:gl511}) see~\cite{knispel}.
\end{remark}

\vspace*{-14pt}

\begin{appendix}\label{app}
\section*{Appendix} Let us finally summarize some technicalities.
\begin{lemma}
\label{th:liptserersatz}
Let $W$ be a two-dimensional Brownian motion on the stochastic base
$(\Omega,\mathcal{G},\mathbb{G},Q)$, and let $\eta$ be a $\mathbb
{G}$-progressively measurable process taking its values in a compact
subset $\Gamma\subset\mathbb{R}^d$.
Moreover, let $(Y_t)_{t\leq T}$ be a continuous process that is a
strong solution of the SDE
%
\begin{equation}
\label{eq:gl83}
dY_t=h(\eta_t,Y_t)\dd t+\sigma\dd W_t,\qquad Y_0=y_0,\qquad \Vert\sigma\Vert>0,
\end{equation}
where the drift function $h\dvtx\Gamma\times\mathbb{R}\rightarrow\mathbb
{R}$ satisfies for all $\eta\in\Gamma$, $y\in\mathbb{R}$
\[
h^2(\eta,y)\leq K^2(1+y^2) \qquad\mbox{for some constant } K.
\]
Then there exists $\delta=\delta(T)>0$ such that $\sup_{t\leq T}E_Q[\exp
(\delta Y_t^2)]<\infty$.
\end{lemma}
\begin{pf} $\!\!\!$The local martingale $B_t\,{:=}\,\Vert\sigma\Vert^{-1}\sigma W_t$,
$t\,{\in}\,[0,T]$, satisfies \mbox{$\langle B\rangle_t\,{=}\,t$}. Thus, $B$ is a
one-dimensional Brownian motion, due to L\'{e}vy's
characterization. In particular, the SDE (\ref{eq:gl83}) can be
rewritten as
%
\begin{equation}
\label{eq:gl393}
dY_t=h(\eta_t,Y_t)\dd t+\Vert\sigma\Vert\dd B_t.
\end{equation}
The proof now follows in two steps: first we argue for a constant
function \mbox{$h(y)\equiv h$}. In that case, the solution to (\ref{eq:gl393})
is given by the Gaussian OU process $
Y_t=e^{ht}(y_0+\int_0^te^{-hs}\Vert\sigma\Vert\dd B_s)$, $t\in[0,T]$,
and the claim follows easily. In a second step, we extend this result
to the general case by a comparison argument. The details are given in
\cite{liptsershiryaev}, Theorem 4.7, restricted to the special
case $h(\eta,y)=h(y)$.
\end{pf}
\begin{lemma}
\label{th:boundedmoments}
Let $(\Omega,\mathcal{G},\mathbb{G},Q)$ be a reference probability
system supporting a two-dimensional Brownian motion $W=(W^1,W^2)$, and
let $\eta$ be a $\mathbb{G}$-progressively measurable process with
values in a compact subset $\Gamma\subset\mathbb{R}^d$. Furthermore, we
suppose that $Y$ is a strong solution to the SDE (\ref{eq:gl83}), where
$h$ is real-valued function such that
\[
\exists K,M>0, \forall\eta\in
\Gamma\dvtx yh(\eta,y)\leq-Ky^2+M,
\]
and where the volatility vector satisfies $\Vert\sigma\Vert>0$. Then it
holds that:
\begin{longlist}
\item There exist constants $C,C_n>0$, $n\in\mathbb{N}$,
such that
\[
\sup_{t\geq0}E_Q[Y_t^{2n}]\leq y_0^{2n}+C_n \quad\mbox{and}\quad \sup
_{t\geq0}E_{Q}[|Y_t|]\leq C(1+|y_0|).
\]
\item For all $k\in\mathbb{R}$, $\sup_{t\geq0}E_Q[\exp
(k\vert Y_t\vert)]<\infty$.
\end{longlist}
In particular, these bounds are uniform among the class of all
progressively measurable $\Gamma$-valued processes $\eta$.
\end{lemma}
\begin{pf}
The proof is rather standard in ergodic control theory and appears in
single components under slight different assumptions in various papers
(see, e.g., \cite{flemingmceneaney} or \cite{hatanagaisheu}). For a
unifying version see \cite{knispel}, Lemma A.2.
\end{pf}
\end{appendix}

\section*{Acknowledgments}
It is a pleasure to thank Hans F\"{o}llmer for helpful discussions. I
am also very grateful for the insightful comments of an anonymous
referee.



%
\printaddresses

\end{document}